\documentclass[11pt,a4paper,twoside]{amsart}
\usepackage{amsmath,amssymb,amsfonts,amsthm,mathrsfs}
\usepackage{times,hyperref,color}

\newcommand{\zero}[1]{\underline{#1}_{\circ}}
\let\mathcal\mathscr
\newtheorem{The}{Theorem}[section]
\newtheorem{Theorem}{Theorem}[section]
\newtheorem{Proposition}[The]{Proposition}
\newtheorem{Lemma}[The]{Lemma}
\newtheorem{Corollary}[The]{Corollary}
\newtheorem{Observation}[The]{Observation}
\theoremstyle{definition}
\newtheorem{Definition}[The]{Definition}
\newtheorem{Remark}[The]{Remark}

\newcommand{\C}{\mathbb{C}}\newcommand{\R}{\mathbb{R}}
\newcommand{\vf}{\vfill\end{document}}

\begin{document}


\title{Cartan equivalences for 
\\
Levi-nondegenerate
hypersurfaces $M^3$ in $\mathbb C^2$ 
\\
belonging to General Class ${\sf I}$}

\author{Masoud Sabzevari}
\address{Department of Pure Mathematics,
University of Shahrekord, 88186-34141 Shahrekord, IRAN and School of
Mathematics, Institute for Research in Fundamental Sciences (IPM), P.
O. Box: 19395-5746, Tehran, IRAN} \email{sabzevari@math.iut.ac.ir}

\author{Jo\"el Merker}
\address{D\'{e}partment de Math\'{e}matiques d'Orsay, B\^{a}timent 425,
Facult\'{e} des Sciences, Universit\'{e} Paris XI - Orsay, F-91405
Orsay Cedex, FRANCE} \email{joel.merker@math.u-psud.fr}

\date{\number\year-\number\month-\number\day}

\maketitle

\begin{abstract}
We develope in great computational details the classical Cartan
equivalence problem for Levi-nondegenerate $\mathcal C^6$-smooth real
hypersurfaces $M^3$ in $\mathbb C^2$, performing all calculations effectively in
terms of a (local) graphing function $\varphi$. In particular, we present
explicitly the unique (complex) essential invariant $\frak J$ of the
problem. Its expansion in terms of the 3-variables function $\varphi$
incorporates millions of differential monomials, while, when $\varphi$ is
assumed to depend only on 2 variables (rigid case), $\frak J$ writes out in
two lines (7 monomials).
\end{abstract}

\pagestyle{headings} \markright{Cartan equivalence problem for $M^3
\subset \mathbb{ C}^2$}

\section{Introduction}

In 1907, Henri Poincar\'{e} \cite{Poincare} initiated the question of
determining whether two given Cauchy-Riemann (CR for short) local
real analytic hypersurfaces in $\mathbb C^2$ can be mapped onto each
other by a certain (local or global) biholomorphism. This problem was
solved later on in 1932 by \'{E}lie Cartan \cite{Cartan} in a
complete way, by importing techniques from his main original impulse
(years 1900--1910) towards general investigations of a large class of
problems which nowadays are known as {\sl Cartan equivalence
problems}, addressing, in many different contexts, equivalences of
submanifolds, of (partial) differential equations, and as well, of
several other geometric structures. Unifying the wide variety of
these seemingly different equivalence problems into a potentially
universal approach, Cartan showed that almost all continuous
classification questions can indeed be reformulated in terms of
specific adapted coframes.

Seeking an equivalence between coframes usually comprises a certain
initial {\sl ambiguity subgroup} $G\subset {\sf Gl}(n)$ related to
the specific features of the geometry under study. The fundamental
general set up is that, for two given coframes $\Omega := \{\omega^1,
\ldots, \omega^n\}$ and $\Omega' := \{{\omega'}^1, \ldots,
{\omega'}^n\}$ on two certain $n$-dimensional manifolds $M$ and $M'$,
there exists a diffeomorphism $\Phi \colon M \longrightarrow M'$
making a geometric equivalence {\em if and only if} there is a
$G$-valued function $g\colon M \rightarrow G$ such that $\Phi^\ast
(\Omega) = g\cdot \Omega'$.

Cartan's `algorithm' (the outcomes of which is often unpredictable)
comprises three interrelated principal aspects: {\sl absorbtion};
{\sl normalization}; {\sl prolongation}.

\smallskip

In brief outline, starting from:
\begin{equation}
\label{equivalence group} \Omega := g\cdot\Omega',
\end{equation}
one has to find the so-called {\sl structure equations} by computing
the exterior differential:
\[
d\Omega = dg\wedge\Omega'+g \cdot d\Omega'.
\]
Inverting~\thetag{\ref{equivalence group}} as $\Omega' =
g^{-1}\,\Omega$, one begins by replacing this in the first term:
\[
dg\wedge\Omega' = dg \wedge g^{-1}\,\Omega = \underbrace{dg\cdot
g^{-1}}_{\text{\sf Maurer-Cartan} \atop \text{\sf matrix}\,\,{\sf
MC}_g} \wedge\, \Omega,
\]
with the standard Maurer-Cartan matrix of the matrix group $G$:
\[
{\sf MC}_g := \bigg( \sum_{k=1}^n\, dg_k^i\, \big(g^{-1}\big)_j^k
\bigg)_{1\leqslant i\leqslant n}^{1\leqslant j\leqslant n} =
\sum_{s=1}^r\, a^i_{js}\,\alpha^s
\]
having $n^2$ entries which express linearly in terms of some basis
$\alpha^1, \dots, \alpha^r$ of left-invariant $1$-forms on $G$, with
$r := \dim_\R\, G$, by means of certain constants $a_{js}^i$. Then
the structure equations become:
\begin{eqnarray*}
d\omega^i = \sum_{j=1}^n\,\sum_{s=1}^r\,a_{js}^i\,\alpha^s \wedge
\omega^j + g\cdot d\Omega' \ \ \ \ \ \ \ \ \ \ \ \ \
{\scriptstyle{(i\,=\,1\,\cdots\,n)}}.
\end{eqnarray*}

Moreover, one has to express the second term $d\Omega'$ above, which
is a $2$-form, as a combination of the $\omega^j \wedge \omega^k$.
Usually, this step is quite costful, computationally speaking. When
one executes this, the appearing (complicated) functions $T_{jk}^i$,
called {\sl torsion coefficients}:
\[
d\omega^i = \sum_{j=1}^n\,\sum_{s=1}^r\,a_{js}^i\,\alpha^s \wedge
\omega^j + \sum_{1\leqslant j<k\leqslant n}\, T_{jk}^i\cdot
\omega^j\wedge\omega^k \ \ \ \ \ \ \ \ \ \ \ \ \
{\scriptstyle{(i\,=\,1\,\cdots\,n)}},
\]
usually reveal {\em appropriate invariants} of the geometric
structure under study.

\smallskip

Then the main thrust of Cartan's approach is that, when one
substitutes each Maurer-Cartan form $\alpha^s$ with
$\alpha^s+\sum_{j=1}^n\,z^s_j\,\omega^j$ for arbitrary
functions-coefficients $z_j^s$, while each torsion coefficient
$T^i_{jk}$ is simultaneously necessarily replaced by $T^i_{jk} +
\sum_{s=1}^r\, \big( a_{js}^i\,z_k^s - a_{ks}^i\,z_j^s \big)$, and
when one does choose the functions-coefficients $z_j^s$ {\em in order
to `absorb' as many as possible torsion coefficients in the
Maurer-Cartan part}, then the remaining, {\sl unabsorbable}, (new,
less numerous) torsion coefficients become {\em true invariants} of
the geometric structure under study. Of course, the `number' of
invariant torsion coefficients is `counted' by means of linear
algebra, usually applying the so-called (non-explicit) Cartan's
Lemma.

Since the remaining torsion coefficients are essential and invariant,
one then {\sl normalizes} them to be equal to a constant, usually
$0$, $1$ or $i$, simply whether or not the group parameters they
contain must be nonzero in the matrix group $G$ to preserve
invertibility. Setting these essential torsions equal to $0$, $1$ or
$i$ then determines some entries of the matrix group $G$, and
therefore decreases the dimension of $G$. In high-level equivalence
problems (\cite{ 5-cubic, Pocchiola}), these {\sl potentially
normalizable} essential torsions are rather numerous and often
overdetermined (unfortunately), hence one is forced to enter more
deeply in explicit computations if one wants to rigorously settle
which group parameters really remain, and which invariants really pop
up. Hopefully at the end of a long procedure, one reduces the
structure group $G$ to dimension $0$, getting a so-called {\sl
$e$-structure}.

But if, as also often occurs, it becomes no longer possible after
several absorption-normalization steps to determine a (reduced) set
of remaining group parameters, then one has to add the rest of
(modified) Maurer-Cartan forms to the initial lifted coframe $\Omega$
and to {\sl prolong} the base manifold $M$ as the product $M^{\sf
pr}:=M\times G$.  Surprisingly, Cartan observed that the solution of
the original equivalence problem can be derived from that of $M^{\sf
pr}$ equipped with the new coframe. Then, one has to restart the
procedure {\it ab initio} with such a new prolonged problem. This
initiates the third essential feature of the equivalence algorithm:
the {\sl prolongation}. For a detailed presentation of Cartan's
method, the reader is referred to \cite{Olver,Gardner,5-cubic}.

Cartan's remarkable achievements were encouraging enough to establish
his elegant geometries, nowadays known as {\sl Cartan geometries}, a
generalization of two seemingly disparate geometries, that of Felix
Klein and that of Bernhard Riemann. For the study of hypersurfaces in
complex Euclidean spaces, Cartan's method was applied later on by
some other mathematicians, {\em e.g.} Chern-Moser \cite{Chern-Moser}
and Tanaka \cite{Tanaka}, but along two seemingly different ways. In
fact, Chern-Moser's work was a fairly direct development of that of
Cartan, while Tanaka's was more algebraically-minded, involving Lie
algebra cohomology, infinitesimal CR automorphisms, and so-called
{\sl Tanaka prolongations}.

\smallskip

Coming to the heart of the matter, let $M^3 \subset \mathbb C^2$ be a
$\mathcal{ C}^6$-smooth Levi-nondegenerate real hypersurface passing
through the origin, in some suitable affine holomorphic coordinates
$(z, w) = (x + i y, u + i\,v)$ represented as the graph of a certain
$\mathcal C^6$-smooth defining function:
\[
v = \varphi(z,\overline{z},u):=z\overline z+{\rm O}(3),
\]
satisfying $\varphi (0)=0$. Our purpose in this paper is to
reformulate Cartan's construction of an $\{ e\}$-structure associated
to such hypersurfaces {\it effectively in terms of the single datum}
$\varphi$ of the problem.

In~\cite{Merker-Sabzevari-CEJM}, inspired by~\cite{ EMS}, we already
performed, within the Tanaka framework, an {\it effective}
construction of a Cartan geometry that is invariantly associated to
such $M^3 \subset \C^2$.  As the main result there, we explicitly
computed the two essential {\it real} curvature coefficients of the
geometry, the vanishing of which characterizes biholomorphic
equivalency of $M$ to the {\sl Heisenberg sphere} $v=z\overline z$
({\it see} Theorem 7.4 in \cite{Merker-Sabzevari-CEJM}). In the
present paper, we have to keep track of how the under consideration
Cartan equivalence problem for real hypersurfaces $M^3$ matches up to
their Cartan-Tanaka geometry. In particular, we will explicitly
observe a close relationship between the single {\it complex}
essential invariant of the equivalence problem and the two {real}
invariants of the Cartan geometry.

\smallskip

As an outline of this paper, first in section \ref{setting-up}, we
set up the equivalence problem for Levi-nondegenerate real
hypersurfaces $M^3\subset\mathbb C^2$ by constructing the necessary
adapted coframe on it. We begin by presenting generators $\mathcal{
L}$ and $\overline{ \mathcal{ L}}$ of $T^{1, 0} M$ and of $T^{0, 1}
M$. Then, the bracket $\mathcal T := i\,\big[ \mathcal
L,\overline{\mathcal L} \big]$ completes a frame on for $\C
\otimes_\R TM$.  Dually, we deduce an initial complex {\em co}frame
$\big\{\rho_0 ,\zeta_0 ,\overline\zeta_0 \big\}$ on $\C \otimes_\R
T^*M$.

Next, we determine the initial ambiguity group for equivalences under
local biholomorphisms:
\begin{eqnarray*}\footnotesize
G:=\left\{g:= \left(\!
\begin{array}{ccc}
{\sf a} & 0 & 0 \\
{\sf b} & {\sf c} & 0 \\
\overline{{\sf b}} & 0 & \overline{\sf c} \\
\end{array}
\!\right), \ \ \ \ {\sf a}\in\mathbb R, \ \ \ {\sf b,c}\in\mathbb
C\right\}.
\end{eqnarray*}

In section \ref{Absorbtion}, we proceed to the equivalence algorithm
by performing the absorbtion-normalization procedure. After
normalizing the group parameter $\sf a$, we continue in section
\ref{Prolongation} by performing a first prolongation.  Namely, we
prolong the equivalence problem of the under consideration
CR-manifolds $M^3$ to that of a certain 7-dimensional prolonged
spaces $M^{\sf pr}:=M^3\times G$ equipped with the initial coframe
$\big\{ \rho_0, \zeta_0, \overline{ \zeta}_0 \big\}$ to which we add
four certain Maurer-Cartan 1-forms
$\alpha,\beta,\overline\alpha,\overline\beta$\,\,---\,\,associated to
certain four remaining group parameters ${\sf b,c},\overline{\sf
b},\overline{\sf c}$\,\,---\,\,and with four new appearing prolonged
group parameters ${\sf r,s},\overline{\sf r},\overline{\sf s}$.
Subsequently, we consider this new prolonged equivalence problem {\it
ab initio}.

The well-known Cartan's Lemma ({\it see} Lemma~\ref{Cartan's Lemma})
also enables us to temporarily bypass some relatively painful
computations ({\it cf.}~Proposition~\ref{Prop-modofoed-MC}), that,
anyway, we do perform later on. After two absorbtions-normalizations
and after one prolongation along the way, the desired equivalence
problem transforms to that of some\,\,---\,\,explicitly
computed\,\,---\,\,eight-dimensional coframe
$\big\{\rho,\zeta,\overline\zeta,\alpha,
\beta,\overline\alpha,\overline\beta,\delta \big\}$ having {\it
e}-structure equations:
\begin{eqnarray*}
\begin{array}{ll}\footnotesize
d\rho=\alpha\wedge\rho+\overline\alpha\wedge\rho+i\,\zeta\wedge\overline\zeta,
&
\\
d\zeta=\beta\wedge\rho+\alpha\wedge\zeta, &
\\
d\overline\zeta=\overline\beta\wedge\rho+\overline\alpha\wedge\overline\zeta,
&
\\
d\alpha=\delta\wedge\rho+2\,i\,\zeta\wedge\overline{\beta}+i\,\overline{\zeta}\wedge\beta,

&
\\
d\beta=\delta\wedge\zeta+\beta\wedge\overline{\alpha}+{\frak
I}\,\overline\zeta\wedge\rho, &
\\
d\overline\alpha=\delta\wedge\rho-2\,i\,\overline\zeta\wedge{\beta}-i\,{\zeta}\wedge\overline\beta,
&
\\
d\overline\beta=\delta\wedge\overline\zeta+\overline\beta\wedge{\alpha}+\overline{\frak
I}\,\zeta\wedge\rho, &
\\
d\delta=\delta\wedge\alpha+\delta\wedge\overline\alpha+i\,\beta\wedge\overline\beta+{\frak
T}\,\rho\wedge\zeta+\overline{\frak T}\,\rho\wedge\overline\zeta, &
\end{array}
\end{eqnarray*}
with the single primary complex invariant:
\[\footnotesize\aligned
{\frak I}&:=-\frac{1}{3}\frac{\overline{\mathcal L}({\mathcal
L}(\overline{\mathcal L}(\overline P)))}{{\sf c}\overline{\sf
c}^3}+\frac{2}{3}\frac{{\mathcal L}(\overline{\mathcal L}(\overline
P))\overline P}{{\sf c}\overline{\sf
c}^3}+\frac{1}{2}\frac{\overline{\mathcal L}(\overline{\mathcal
L}({\mathcal L}(\overline P)))}{{\sf c}\overline{\sf
c}^3}-\frac{7}{6}\frac{\overline{\mathcal L}({\mathcal L}(\overline
P))\overline P}{{\sf c}\overline{\sf c}^3}-
\\
&\ \ \ \ \ -\frac{1}{6}\frac{{\mathcal L}(\overline
P)\overline{\mathcal L}(\overline P)}{{\sf c}\overline{\sf
c}^3}+\frac{1}{3}\frac{{\mathcal L}(\overline P)\overline P^2}{{\sf
c}\overline{\sf c}^3},
\endaligned
\]
in which the fundamental function $P$ can expresses explicitly in
terms of the single datum $\varphi$ of the problem as:
\[
P:=\frac{\ell_z-\ell A_u+A\ell_u}{\ell},
\]
where:
\[
A := \frac{i\,\varphi_z}{1-i\,\varphi_u} \ \ \ \ \ \ \ \ \ \ \ \
\text{\rm and where:} \ \ \ \ \ \ \ \ \ \ \ \ \ell := i\,\big(
\overline{A}_z + A\,\overline{A}_u - A_{\overline{z}} -
\overline{A}\,A_u \big),
\]
this last {\sl Levi factor} $\ell$ being nowhere vanishing, because
we assume $M$ to be Levi nondegenerate. Furthermore, the other
secondary invariant $\frak T$ can be expressed in terms of the first
one $\frak I$ as:
\[
{\frak T}=\frac{1}{\overline{\sf c}}\,\bigg(\overline{\mathcal
L}(\overline{\frak I})-\overline P\,\overline{\frak
I}\bigg)-i\,\frac{\sf b}{{\sf c}\overline{\sf c}}\overline{\frak I}.
\]

Finally in section \ref{Cartan geometry}, we turn to a brief
discussion of the Cartan-Tanaka geometry of the under consideration
hypersurfaces $M^3$ and\,\,---\,\,being aware of the results of the
papers \cite{EMS, Merker-Sabzevari-CEJM}\,\,---\,\,we observe that
the equivalence problem matches up to their Cartan geometry so that
the complex essential primary invariant $\frak J$ can be reexpressed
effectively in terms of the two (real) essential primary invariants
we obtained there (this also matches up with the results of~\cite{
Merker}).

\begin{Theorem} (see Theorem
\ref{Th-Geometry} at the end) For Levi-nondegenerate $\mathcal
C^6$-smooth real hypersurfaces $M^3\subset\mathbb C^2$, the following
relation holds between the essential complex invariant $\frak J$ of
their equivalence problem and the essential real invariants
$\mathbf{\Delta}_1$ and $\mathbf{\Delta}_2$ of their Cartan geometry:
\[
\frak I=\frac{4}{{\sf c}\overline{\sf
c}^3}\big(\mathbf{\Delta}_1+i\,\mathbf{\Delta}_4\big).
\]
\end{Theorem}

We close up this introduction by mentioning that, although it is well
known that a close relationship exists between equivalences of
hypersurfaces $M^3 \subset \C^2$ and second-order ordinary
differential equations (\cite{ Cartan, Chern-Moser, Jacobowitz, GTW,
Sparling-Nurowski, Merker}), and although the (nonexplicit) geometric
features of the results we present here are well known too (but often
with hidden computations), a completely effective and systematic
presentation of the related (complicated) computational aspects is
necessary to understand in a deeper way the core of Cartan's method.

\smallskip

In fact, the present (preliminary) paper was written up in order to
serve as a ground-companion to much higher level explorations of
equivalence problems for embedded CR structures, that will appear
soon (\cite{ 5-cubic, Pocchiola}). Intentionally, we endeavour here
to develope our systematic computational formalism at first for the
simplest known CR structures $M^3 \subset \C^2$, before applying it
to more delicate $5$-dimensional real analytic CR structures.

\smallskip

The remarkable works of Beloshapka~\cite{Beloshapka1997,
Beloshapka2002, Beloshapka2004, Beloshapka2007, BES} have shown that
there exists a wealth of model CR-generic submanifolds whose algebras
of infinitesimal CR automorphisms have been computed {\em explicitly}
there, and this paper together with~\cite{ BES, 5-cubic, Pocchiola}
are a very first step in the Cartan-like study of the
geometry-preserving deformations of just a few of these models, with
a door potentially open towards the exploration of a great number of
higher models with a similar emphasis on {\em effectiveness}.

\section{Setting up the equivalence problem}
\label{setting-up}

Our aim in this section is to construct\,\,---\,\,in terms of a
certain fundamental graphing function $\varphi$\,\,---\,\,an initial
complex coframe on the under consideration three dimensional
CR-manifold $M^3\subset\mathbb C^2$, and next to set up the related
equivalence problem. First, let us consider this approach dually,
namely by constructing a local {\em frame} on $M^3$.

\subsection{Local frame adapted to $3$-dimensional
embedded CR structures} Consider therefore a local $\mathcal{
C}^6$-smooth hypersurface $M^3 \subset \mathbb C^2$ passing through
the origin. In some suitable affine holomorphic coordinates $(z, w) =
(x + i y, u + i\,v)$ adapted so that $T_0 M^3 = \{ v = 0\}$, the
implicit function theorem enables one to represent $M^3$ as a graph
over the $(x, y, u)$-space. Since any function of $(x, y, u) = \big(
\frac{ z + \overline{ z}}{ 2}, \, \frac{ z + \overline{ z}}{ 2 i}, \,
u)$ can be considered as one of $( z, \overline{ z}, u)$, the graph in
question may be thought of as being of the form:
\[
v = \varphi(z,\overline{z},u),
\]
for some $\mathcal{ C}^6$ function $\varphi$ satisfying $\varphi (0) =
\varphi_z ( 0) = \varphi_{ \overline{ z}} ( 0) = \varphi_u ( 0)$. In
the sequel, all appearing invariant objects\,\,---\,\,vector fields,
differential forms, torsion coefficients, essential
functions\,\,---\,\,will depend only on $\varphi$ and its partial
derivatives with respect to the three (complex and real) initial
coordinates $(z, \overline{ z}, u)$, the latter being understood as
{\sl intrinsic} coordinates on $M^3$.

According to \cite{ Jacobowitz, Merker-Sabzevari-CEJM}, a local $(1,
0)$ vector field on $\C^2$ defined near the origin:
\[
\mathcal{L} := \frac{\partial}{\partial z} + {\tt
A}\,\frac{\partial}{\partial w}
\]
is {\sl tangent} to $M^3$ if and only if, on restriction to $M^3$,
its coefficient ${\tt A}$ satisfies:
\[
\aligned 0 & = \mathcal{L} \Big( -
{\textstyle{\frac{w-\overline{w}}{2i}}} + \varphi\big(z,\overline{z},
{\textstyle{\frac{w+\overline{w}}{2}}}\big) \Big)
\\
& = -\,\frac{1}{2i}\,{\tt A} + \frac{1}{2}\,{\tt A}\,\varphi_u +
\varphi_z.
\endaligned
\]
For this to hold true, it suffices to set:
\[
{\tt A} := \frac{-2\,\,\varphi_z}{i+\,\varphi_u},
\]
which is thus {\em de facto} a function of only $(z, \overline{ z},
u)$. Furthermore, restricting $\mathcal L$ to $M^3$, one must simply
and only drop the (extrinsic) vector field $\frac{ \partial}{ \partial v}$:
\[
\aligned \mathcal{L}\big\vert_M & = \frac{\partial}{\partial z} +
{\tt A}\, \bigg( \frac{1}{2}\, \frac{\partial}{\partial u} -
\zero{\frac{i}{2}\, \frac{\partial}{\partial v}} \bigg)
\\
& = \frac{\partial}{\partial z} - \frac{\,\varphi_z}{i+\varphi_u}\,
\frac{\partial}{\partial u}.
\endaligned
\]

Now, it will be convenient to introduce an extra notation for the
appearing coefficient of $\frac{ \partial}{ \partial u}$, say:
\begin{equation}
\label{A} A := \frac{i\,\varphi_z}{1-i\,\varphi_u},
\end{equation}
not to be confused with ${\tt A} = 2\, A$, which, anyway, will be
left aside from now on.

Thus intrinsically on $M^3$, the CR-structure induced by the ambient
$\C^2$ on $M^3$ is encoded by the complex $(1, 0)$ vector field $\mathcal{ L}$
and its conjugate $\overline{ \mathcal{ L}}$:
\[
\mathcal{L} = \frac{\partial}{\partial z} +
A\,\frac{\partial}{\partial u} \ \ \ \ \ \ \ \ \ \ \text{\rm and} \ \
\ \ \ \ \ \ \ \ \overline{\mathcal{L}} =
\frac{\partial}{\partial\overline{z}} + \overline{A}\,
\frac{\partial}{\partial u}.
\]
In this set up, the non-vanishing property of the Lie bracket:
\[
\big[\mathcal{L},\,\overline{\mathcal{L}}\big] = \big( \overline{A}_z
+ A\,\overline{A}_u - A_{\overline{z}} - \overline{A}\,A_u \big)\,
\frac{\partial}{\partial u}
\]
at any point of $M^3$ indicates precisely that $M^3$ is {\sl Levi
nondegenerate} at every point, an assumption
that will be held throughout. Since it is slightly
better\,\,---\,\,for convenience reasons\,\,---\,\,to deal with {\em
real} functions, we introduce the fundamental {\sl Levi
factor}:
\begin{equation}
\label{ell} \ell := i\,\big( \overline{A}_z + A\,\overline{A}_u -
A_{\overline{z}} - \overline{A}\,A_u \big),
\end{equation}
so that the reality of $\ell\, \frac{\partial}{ \partial u}$ in the
first structural Lie bracket relation, viewed again in this
abbreviated way $[\mathcal{L},\,\overline{\mathcal{L}}] =
-i\,\ell\,\frac{\partial}{\partial u}$, shows now well that the $-i$
mere factor on the right provides the pure imaginarity of the bracket
in question:
\[
\overline{\big[\mathcal{L},\overline{\mathcal{L}}\big]}
=
-\,\big[\mathcal{L},\,\overline{\mathcal{L}}\big].
\]

For normalization reasons,
it is furthermore natural to introduce
the auxiliary {\it real} field:
\[
\mathcal{T} := \ell\,\frac{\partial}{\partial u},
\]
which is the suitable multiple of $\frac{ \partial}{ \partial u}$ insuring
that the bracket:
\[
\big[ \mathcal{L},\,\overline{\mathcal{L}}\big] = -i\,\mathcal{T}
\]
makes the coefficient-function in front of $\mathcal{ T}$ to
become a plain constant.

Now, in terms of what will be called the {\sl complex initial frame}
on $M^3$ (written in the following order):
\[\footnotesize
\boxed{ \aligned \mathcal{T} & :=
i\,\big(\overline{A}_z+A\,\overline{A}_u - A_{\overline{z}} -
\overline{A}\,A_u\big)\, \frac{\partial}{\partial u},
\\
\mathcal{L} & := \frac{\partial}{\partial z} +
A\,\frac{\partial}{\partial u},
\\
\overline{\mathcal{L}} & := \frac{\partial}{\partial\overline{z}} +
\overline{A}\, \frac{\partial}{\partial u},
\endaligned
}
\]
it remains to also take up the two remaining\,\,---\,\,yet
uncomputed\,\,---\,\,brackets.

Simple computations show that we have:
\[
\big[\mathcal{T},\,\mathcal{L}\big] = -\,P\,\mathcal{T} \ \ \ \ \ \ \
\ \ \text{\rm and} \ \ \ \ \ \ \ \ \
\big[\mathcal{T},\,\overline{\mathcal{L}}\big] =
-\,\overline{P}\,\mathcal{T},
\]
for a certain (universal) rational function $P$ of the second-order
jet $J_{ z, \overline{ z}, u} \big( A, \overline{ A}\big)$
given by:
\[
P:=\frac{\ell_z-\ell A_u+A\ell_u}{\ell}.
\]

This function $P$ could be completely
expanded in terms of the graphing function
$\varphi$, for in the notation
of~\cite{Merker-Sabzevari-CEJM}, one checks that:
\[
P
=
{\textstyle{\frac{1}{2}}}\,\Phi_1
-
{\textstyle{\frac{i}{2}}}\,\Phi_2,
\]
with the full, one-page long, expansions of
(the numeratorof) $\Phi_1$ and
$\Phi_2$ in terms
of $J_{x, y, u}^3 \varphi$ being provided on page~42 of
the extensive {\small\sf arxiv.org} version
of~\cite{Merker-Sabzevari-CEJM}.
Because the computations unavoidably {\sl explode}
when one performs them in terms of $\varphi$
({\em cf.} the end of~\cite{Merker-Sabzevari-CEJM}),
it is advisable to reset oneself at the level
of just $P$, aiming nevertheless to perform everything
which will follow in terms of $P$, granted that $P$
is explicit with respect to $\varphi$.

Notice {\em passim} that the above two structural bracket relations
are conjugate to each other, just because $\overline{ \mathcal{ T}}
={ \mathcal{ T}}$. Furthermore:

\begin{Lemma}
\label{LP} One has the reality condition:
\[
\mathcal{L}(\overline P)
=
\overline{\mathcal L}(P).
\]
\end{Lemma}

\proof
The already presented expressions simply give:
\[
\footnotesize\aligned \big[\overline{\mathcal
L},\underbrace{[\mathcal T,\mathcal L]}_{-P\mathcal
T}\big]&=-\overline{\mathcal L}(P)\mathcal T-P\overline P\mathcal T,
\\
\big[\mathcal L,\underbrace{[\overline{\mathcal L},\mathcal
T]}_{\overline P\mathcal T}\big]&=\mathcal L(\overline P)\mathcal
T+\overline PP,
\endaligned
\]
and thanks to the Jacobi identity, one obtains:
\[
\small \aligned -\overline{\mathcal L}(P)\mathcal T+\mathcal
L(\overline P)\mathcal T & = \big[\overline{\mathcal L},[\mathcal
T,\mathcal L]\big]+\big[\mathcal L,[\overline{\mathcal L},\mathcal
T]\big] = -\big[\mathcal T,\underbrace{[\mathcal L,\overline{\mathcal
L}]}_{\mathcal T}\big]=0,
\endaligned
\]
which visibly yields the desired equality $\mathcal
L(\overline P)=\overline{\mathcal L}(P)$.
\endproof

\subsection{Setting up of an initial Cartan coframe}

All these preliminary normalizations were done in advance to fit
dually with a pleasant collection of 1-forms. Indeed, on the natural
agreement that the coframe $\{du,\,dz,\,d\overline{z}\}$ is dual to
the frame $\{ {\textstyle{\frac{\partial}{\partial u}}},\,
{\textstyle{\frac{\partial}{\partial z}}},\,
{\textstyle{\frac{\partial}{\partial\overline{z}}}} \}$, let us
introduce the coframe:
\[
\big\{\rho_0,\,\zeta_0,\,\overline{\zeta}_0\big\} \ \ \ \text{\rm
which is dual to the frame} \ \ \
\big\{\mathcal{T},\,\mathcal{L},\,\overline{\mathcal{L}}\big\}.
\]
that is to say which satisfies by definition:
\[
\begin{array}{ccc}
\rho_0(\mathcal{T})=1
\ \ \ & \ \ \
\rho_0(\mathcal{L})=0
\ \ \ & \ \ \
\rho_0(\overline{\mathcal{L}})=0,
\\
\zeta_0(\mathcal{T})=0
\ \ \ & \ \ \
\zeta_0(\mathcal{L})=1
\ \ \ & \ \ \
\zeta_0\big(\overline{\mathcal{L}}\big)=0,
\\
\overline{\zeta_0}(\mathcal{T})=0
\ \ \ & \ \ \
\overline{\zeta_0}(\mathcal{L})=0
\ \ \ & \ \ \
\overline{\zeta_0}\big(\overline{\mathcal{L}}\big)=1.
\end{array}
\]

Using the above expressions of our three vector fields $\mathcal{
T}$, $\mathcal{ L}$, $\overline{ \mathcal{ L}}$,
we see that the
three dual $1$-forms have the following simple explicit expressions in
terms of the function $A$\,\,---\,\,strictly speaking in terms of the
defining function $\varphi$\,\,---\,\,:
\begin{equation}
\label{rho0} \rho_0 :=
\frac{du-A\,dz-\overline{A}\,d\overline{z}}{\ell},\ \ \ \ \ \ \ \
\zeta_0 := dz,\ \ \ \ \ \ \overline{\zeta}_0 := d\overline{z}.
\end{equation}
In order to find the exterior differentiations of these initial
1-forms, an application of the so-called Cartan formula $d\omega ( \mathcal{
X}, \mathcal{ Y}) = \mathcal{ X} \big( \omega ( \mathcal{ Y})\big) -
\mathcal{ Y} \big( \omega ( \mathcal{ X})\big) - \omega\big( [
\mathcal{ X}, \mathcal{ Y}]\big)$ implies that:

\begin{Lemma}
Given a frame $\big\{ \mathcal{ L}_1, \dots, \mathcal{ L}_n\big\}$ on
an open subset of $\R^n$ enjoying the Lie structure:
\[
\big[\mathcal{L}_{i_1},\,\mathcal{L}_{i_2}\big] =
\sum_{k=1}^n\,a_{i_1,i_2}^k\,\mathcal{L}_k \ \ \ \ \ \ \ \ \ \ \ \ \
{\scriptstyle{(1\,\leqslant\,i_1\,<\,i_2\,\leqslant\,n)}},
\]
where the $a_{ i_1, i_2}^k$ are functions on $\R^n$, the dual coframe
$\{ \omega^1, \dots, \omega^n \}$ satisfying by definition $\omega^k
\big( \mathcal{ L}_i\big) = \delta_i^k$ enjoys a quite similar
Darboux-Cartan structure, up to an overall minus sign:
\[
d\omega^k = - \sum_{1\leqslant i_1<i_2\leqslant n}\,
a_{i_1,i_2}^k\,\omega^{i_1}\wedge\omega^{i_2} \ \ \ \ \ \ \ \ \ \ \ \
\ {\scriptstyle{(k\,=\,1\,\cdots\,n)}}.
\]
\end{Lemma}

To apply this lemma, it is convenient to consider the auxiliary array:

\[\footnotesize
\begin{array}{ccccccccccc}
& & \mathcal{T} & & \overline{\mathcal{L}} & & \mathcal{L} & &
\\
& & \boxed{d\rho_0} & & \boxed{d\overline{\zeta_0}} & &
\boxed{d\zeta_0} & &
\\
\big[\mathcal{T},\,\overline{\mathcal{L}}\big] & = &
-\,\overline{P}\cdot\mathcal{T} & + & 0 & + & 0 &
\boxed{\rho_0\wedge\overline{\zeta_0}}
\\
\big[\mathcal{T},\,\mathcal{L}\big] & = & -\,P\cdot\mathcal{T} & + &
0 & + & 0 & \boxed{\rho_0\wedge\zeta_0}
\\
\big[\overline{\mathcal{L}},\,\mathcal{L}\big] & = &
+i\,\cdot\mathcal{T} & + & 0 & + & 0 & \
\boxed{\overline{\zeta_0}\wedge\zeta_0}\,,
\end{array}
\]
in which, by reading the three columns, we deduce visually the
initial Darboux-Cartan structure in terms of our
basic, single function $P$:
\begin{equation}
\label{d0} \boxed{ \aligned d\rho_0 &
=
P\,\rho_0\wedge\zeta_0
+
\overline{P}\,\rho_0\wedge\overline{\zeta_0}
+
i\,\zeta_0\wedge\overline{\zeta}_0,
\\
d\zeta_0 & = 0,
\\
d\overline{\zeta_0} & = 0.
\endaligned}
\end{equation}

\subsection{Complex structure
on the kernel of the contact $1$-form $\rho_0$}

We end up this preparative part by a thoughtful summary which will
offer the natural geometric meaning of $\rho_0$. The defining
equation of $M^3$ may be understood as:
\[
r = 0 \ \ \ \ \ \ \text{\rm with}\ \ \ \ r=r(z,\overline{z},u,v) :=
-v+\varphi(z,\overline{z},u).
\]
Given any function $G = G ( z, \overline{ z}, w, \overline{ w})$, one
classically defines its $(1, 0)$ and $(0, 1)$ differentials
respectively by:
\[
\partial G
:= G_z\,dz+G_w\,dw \ \ \ \ \ \ \ \ \ \text{\rm and} \ \ \ \ \ \ \ \ \
\overline{\partial}G := G_{\overline{z}}\,d\overline{z} +
G_{\overline{w}}\,d\overline{w},
\]
and one easily checks that its complete real differential:
\[
dG = G_x\,dx+G_y\,dy+G_u\,du+G_v\,dv
\]
is the plain sum of these two holomorphic and antiholomorphic
differentials:
\[
dG =
\partial G
+ \overline{\partial}G.
\]

\begin{Lemma}
With $r = 0$ being {\em any} real defining equation for a $\mathcal{
C}^1$ hypersurface $M^3 \subset \C^2$, the restriction to $M^3$ of
the $(1, 0)$ form $i\, \partial r$, namely:
\[
\varrho := i\,\partial r\big\vert_M
\]
is a real form on $M^3$:
\[
\varrho
=
\overline{\varrho}.
\]
Moreover, at every point $p \in M$, the real kernel of $\varrho$ in
$T_pM$ identifies with the complex tangent bundle at $p$:
\[
\big\{ X_p\in\,T_pM\colon\, \varrho(X_p)=0 \big\} = T_p^cM,
\]
while its kernel in the complexified tangent bundle $\C \otimes T_pM$
identifies with $\C \otimes T_p^c M$:
\[
\big\{ \mathcal{X}_p\in\,\C\otimes T_pM\colon\,
\varrho(\mathcal{X}_p)=0 \big\} = \C\otimes T_p^cM = T_p^{1,0}M
\oplus T_p^{0,1}M.
\]
\end{Lemma}

\proof For the first part of the assertion, since $r|_{M^3} \equiv
0$, then on restriction to $M^3$ we also have $dr = 0$ which means
$\partial r = - \overline{
\partial} r$. Hence the $i$ factor in $\varrho$ in front of $\partial
r$ makes it real. For the rest, {\it see} ~\cite{Jacobowitz}, page
25.
\endproof

To go into this lemma in detail, with
$r(z,\overline{z},u,v)=-v+\varphi(z,\overline{z},u)$ and with
$w=u+iv$, we have:
\[
\aligned
dw&=du+i\,dv=du+i\,d\varphi(z,\overline{z},u)=du+i\big(\varphi_z\,dz+\varphi_{\overline{z}}\,d\overline{z}+\varphi_u\,du\big),
\endaligned
\]
and hence the expression of $\varrho$ can be expressed in terms of
the functions $\varphi$:
\begin{equation}
\small\label{varrho} \aligned \ \ \ \ \ \ \ \ \ \ \ \ \ \ \
\varrho&=i\,\partial r|_{M^3}=i\big(r_z\,dz+r_w\,dw\big)|_{M^3}
\\
&=i\big(\varphi_z\,dz+({\textstyle{\frac{1}{2}}}\,\varphi_u+{\textstyle{\frac{i}{2}}})\,dw\big)
\\
&=i\big(\varphi_z\,dz+({\textstyle{\frac{1}{2}}}\,\varphi_u+{\textstyle{\frac{i}{2}}})(du+i\,\varphi_z\,dz+i\,\varphi_{\overline
z}\,d\overline{z}+i\,\varphi_u\,du)\big)
\\
&=\big( {-\textstyle{\frac{1}{2}}} -
{\textstyle{\frac{1}{2}}}\,(\varphi_u)^2 \big)\,du + \big(
{\textstyle{\frac{i}{2}}}\,\varphi_z -
{\textstyle{\frac{1}{2}}}\,\varphi_z\,\varphi_u \big)\,dz + \big(
-{\textstyle{\frac{i}{2}}}\,\varphi_{\overline{z}} -
{\textstyle{\frac{1}{2}}}\,\varphi_{\overline{z}}\,\varphi_u
\big)\,d\overline{z}.
\endaligned
\end{equation}
Furthermore, a plain computations show that ({\it see}
\thetag{\ref{A}}, \thetag{\ref{ell}} and \thetag{\ref{rho0}} for the
expressions):
\begin{equation}
\label{rho-varrho}
\rho_0=-\frac{1}{\ell}\,\frac{2}{1+\varphi_u^2}\,\varrho.
\end{equation}
Then, non-vanishing property of the Levi factor $\ell$ also implies
the equality:
\[
{\rm Ker}(\varrho)
=
{\rm Ker}(\rho_0).
\]

\subsection{Differential facts about CR equivalences}

Now, we explain how one may launch Cartan's method in the
case under study, namely for deformations of the Heisenberg sphere:
\begin{equation}
\label{Heisenberg-sphere} w-\overline{w} = 2i\,z\overline{z},
\end{equation}
that are {\sl geometry-preserving} in the sense
that Levi nondegeneracy is preserved.

Consider therefore two Levi-nondegenerate real hypersurfaces of class
$\mathcal{ C}^6$, represented in two systems of coordinates $(z, w)$
and $(z', w')$ as graphs:
\[
M^3\colon\ \ \ 0 = -v+\varphi(z,\overline{z},u) \ \ \ \ \ \ \
\text{\rm and} \ \ \ \ \ \ \ M'^3\colon\ \ \ 0 =
-v'+\varphi'(z',\overline{z}',u'),
\]
for two certain functions, normalized in advance so that $\varphi :=
z\overline{ z} + {\rm O} ( 3)$ and $\varphi' := z'\overline{ z}' +
{\rm O} ( 3)$. The general problem is to discover {\it when}, and if
so {\it how}, the two CR hypersurfaces are equivalent through a local
ambient biholomorphic map:
\[
(z,w) \longmapsto (z',w') = \big(z'(z,w),\,w'(z,w)\big)
\]
of $\C^2$. This is nothing else than saying that such a map should
send any point of $M^3$ to some determinate point of $M'^3$. In other
words, one should have $v' = \varphi'( z', \overline{ z}', u')$ as
soon as $v = \varphi ( z, \overline{ z}, u)$.

Then a well known
simple fact (Lemma 1.2.3 page 47 of \cite{Stormark})
insures that $M$ is sent to $M'$
if and only if there exists a {\it real-valued}
function $a = a ( z, w)$ defined in a neighborhood of the origin in
$\C^2$ so that:
\[
-\,v'+\varphi'(z',\overline{z}',u')
\big\vert_{(z',w')=(z'(z,w),w'(z,w))} \equiv a(z,,\overline{z},w,
\overline{w})\cdot
\big(-v+\varphi(z,\overline{z},u)\big),
\]
identically as functions of the {\em four} real coordinates of
$\mathbb C^2$. For easier reading, we shall drop the mention of this
{\it pullback} and simply write down:
\[
-v'+\varphi'(z',\overline{z}',u') =
a\,\big(-v+\varphi(z,\overline{z},u)\big),
\]
or even in a shorter way: $r' = a\, r$. We now clearly see that $r = 0$
implies $r' = 0$, namely that points of $M^3$ are sent to points of
$M'^3$. But now, the two fundamental $1$-forms $\varrho=i\,
\partial r \big\vert_M$ and $\varrho'=i\, \partial r' \big\vert_{ M'}$ in
the two spaces happen to be {\em real multiples of each other}:
\[
i\,\partial r'\big\vert_{M'} = a\,i\,\partial r\big\vert_M +
\zero{r\,i\,\partial a\big\vert_M},
\]
through the same function $a$.

Of course such a function $a$ highly
depends on the equivalence $(z, w) \to (z', w')$ between $M^3$ and
$M'^3$, when it exists, but the idea of Cartan is to consider it as
some {\em unknown}. Taking the relationship \thetag{\ref{rho-varrho}} into
account, the already obtained equality $\varrho'=a\,\varrho$ can be
slightly adjusted (with same notation for a new
function $a$) into the form:
\[
\rho':=a\cdot\rho
\]
for some unknown real-valued function $a:=a(z,\overline z,u)$.

\subsection{Associated ambiguity matrix}

Next, let us construct the associated {\sl ambiguity matrix} which
encodes holomorphic equivalence of two hypersurfaces $M^3$ and
$M'^3$, recently equipped with two coframes:
\[
\big\{\rho_0,dz,d\overline{ z}\big\} \ \ \ \ \ \ \ \ \ \text{\rm and}
\ \ \ \ \ \ \ \ \ \big\{\rho_0',dz',d\overline{z'}\big\}.
\]

On restriction to $M^3$, we have:
\[
z' = z'\big(z,\,u+i\,\varphi(z,\overline{z},u)\big),
\]
whence differentiation using the general formula $dg = g_z dz + g_{
\overline{ z}}\, d\overline{ z} + g_u\, du$ gives ({\it see}
\thetag{\ref{varrho}}):
\begin{equation}
\footnotesize \label{dz'} \aligned dz' & =
\big(z_z'+i\,z_w'\,\varphi_z\big)\,dz +
\big(i\,z_w'\,\varphi_{\overline{z}}\big)\,d\overline{z} +
\big(i\,z_w'\,\varphi_u+z_w'\big)\,du
\\
& = \big(z_z'+i\,z_w'\,\varphi_z\big)\,dz + z_w' \big\{\underline{
i\,\varphi_{\overline{z}}\,d\overline{z} + (i\,\varphi_u+1)\,du
}\big\}.
\endaligned
\end{equation}
On the other hand, multiplying by some (innocuous) complex multiple
the fundamental $1$-form $\varrho = i \, \partial r \big\vert_M$, we
also have:
\[\aligned
\frac{-2\,(1+i\,\varphi_u)}{1+(\varphi_u)^2}\, \varrho =
(1+i\,\varphi_u)\,du + i\, \varphi_{\overline{z}}\,d\overline{z} +
\frac{\varphi_z(-i+\varphi_u)}{1-i\,\varphi_u}\,dz,
\endaligned
\]
which enables us to substitute the (underlined) 1-form that we left
in braces after $z_w'$ just above (we also replace $\varrho =
-\frac{1}{2}\, \ell(1+ (\varphi_u)^2)\, \rho_0$ in terms of $\rho_0$,
{\it see} \thetag{\ref{rho-varrho}}) as:
\[\footnotesize
\aligned i\,\varphi_{\overline{z}}\,d\overline{z} +
(i\,\varphi_u+1)\,du & =
-\frac{2\,(1+i\,\varphi_u)}{1+(\varphi_u)^2}\,\varrho -
\frac{\varphi_z(-i+\varphi_u)}{1-i\,\varphi_u}\,dz  =
\ell\,(1+i\,\varphi_u)\,\rho_0 -
\frac{\varphi_z(-i+\varphi_u)}{1-i\,\varphi_u}\,dz.
\endaligned
\]
This implies from \thetag{\ref{dz'}} that $dz'$ is a linear
combination\,\,---\,\,with some complicated coefficients\,\,---\,\,of
$dz$ and of $\rho$, {\em without $d\overline{ z}$ component}:
\[\footnotesize
dz' = \bigg(
\underbrace{z_w'\,\ell\,(1+i\,\varphi_u))}_{=:b(z,\overline{z},v)}
\bigg)\,\rho_0 + \bigg( \underbrace{ z_z'+i\,z_w'\,\varphi_z -z_w'\,
\frac{\varphi_z(-i+\varphi_u)}{1-i\,\varphi_u}}_{=:c(z,\overline{z},v)}\bigg)\,dz.
\]
We thus have obtained:

\begin{Proposition}
Two local $\mathcal{ C}^1$
real hypersurfaces $M^3$ and $M'^3$ of $\mathbb C^2$ are
equivalent through some biholomorphism whenever their two
corresponding fundamental coframes:
\[
\big\{ \rho_0,\,\zeta_0 = dz,\, \overline{\zeta}_0=d\overline{z}_0
\big\} \ \ \ \ \ \ \ \ \ \ \text{\rm and} \ \ \ \ \ \ \ \ \ \ \big\{
\rho_0',\,\zeta_0' = dz',\, \overline{\zeta_0'}=d\overline{z}_0'
\big\}
\]
are mapped one to another by means of a certain matrix of functions:
\[
\left(\!\!
\begin{array}{c}
\rho_0'
\\
\zeta_0'
\\
\overline{\zeta_0'}
\end{array}
\!\!\right) = \left(\!\!
\begin{array}{ccc}
a & 0 & 0
\\
b & c & 0
\\
\overline{b} & 0 & \overline{c}
\end{array}
\!\!\right)\, \left(\!\!
\begin{array}{c}
\rho_0
\\
\zeta_0
\\
\overline{\zeta}_0
\end{array}
\!\!\right),
\]
in which $a := a( z, \overline{ z}, v)$ is a real-valued function on
$M^3$, and where $b := b ( z, \overline{ z}, v)$ and $c := c ( z,
\overline{ z}, v)$ are both complex-valued.
\qed
\end{Proposition}

\subsection{The related structure group}

As we saw, when a CR equivalence exists, the functions $a$, $b$ and
$c$ depend\,\,---\,\,in a somewhat complicated way\,\,---\,\,upon the CR
equivalence, whose existence is under question! The gist of Cartan's
method is to consider these functions {\em as new unknowns}, hence to
add them as extra {\em group variables}. So we consider the subgroup
of matrices inside ${\sf GL}_3 ( \C)$:
\[
\left(\!\!
\begin{array}{ccc}
{\sf a} & 0 & 0
\\
{\sf b} & {\sf c} & 0
\\
\overline{\sf b} & 0 & \overline{\sf c}
\end{array}
\!\!\right),
\]
where now ${\sf a} \in \R$, ${\sf b} \in \C$, ${\sf c} \in \C$ are
arbitrary parameters and we consider the so-called {\sl lifted
coframe} on the eight-dimensional space $(z, \overline{ z}, u, {\sf
a}, {\sf b}, \overline{\sf b}, {\sf c}, \overline{\sf c})$:
\[
\left(\!\!
\begin{array}{c}
\rho
\\
\zeta
\\
\overline{\zeta}
\end{array}
\!\!\right) := \left(\!\!
\begin{array}{ccc}
{\sf a} & 0 & 0
\\
{\sf b} & {\sf c} & 0
\\
\overline{\sf b} & 0 & \overline{\sf c}
\end{array}
\!\!\right) \left(\!\!
\begin{array}{c}
\rho_0
\\
\zeta_0
\\
\overline{\zeta}_0
\end{array}
\!\!\right),
\]
that is to say:
\[
\aligned \rho & = {\sf a}\,\rho_0,
\\
\zeta & = {\sf b}\,\rho_0 + {\sf c}\,\zeta_0,
\\
\overline{\zeta} & = \overline{\sf b}\,\rho_0 + \overline{\sf
c}\,\overline{\zeta}_0.
\endaligned
\]
Of course, the $1$-form $\rho$ is real and the $\overline{ \zeta}$ is
the conjugate of $\zeta$.

So far, we have provided the necessary data for launching the Cartan
algorithm of equivalence. Next, we have
to perform {\sl
normalization}, {\sl absorption} and {\sl prolongation}.

\section{Absorbtion and Normalization}
\label{Absorbtion}

Associated to the equivalence problem for real hypersurface
$M^3\subset\mathbb C^2$, we set up the structure matrix group:
\begin{eqnarray*}\footnotesize
G:=\left\{g:=
\left(%
\begin{array}{ccc}
{\sf a} & 0 & 0 \\
{\sf b} & {\sf c} & 0 \\
\overline{{\sf b}} & 0 & \overline{\sf c} \\
\end{array}%
\right), \ \ \ \ {\sf a}\in\mathbb R, \ \ \ {\sf b,c}\in\mathbb
C\right\}.
\end{eqnarray*}
The lifted coframe writes out as:
\begin{eqnarray}\footnotesize
\label{lifted-coframe}
\left(%
\begin{array}{c}
\rho \\
\zeta \\
\overline{\zeta} \\
\end{array}%
\right) :=g.\left(%
\begin{array}{c}
\rho_0 \\
\zeta_0 \\
\overline{\zeta}_0 \\
\end{array}%
\right)
=\left(%
\begin{array}{c}
{\sf a}\rho_0 \\
{\sf b}\rho_0+{\sf c}\zeta_0 \\
\overline{{\sf b}}\rho_0+\overline{{\sf c}}\overline{\zeta}_0 \\
\end{array}%
\right).
\end{eqnarray}

Applying the differential operator $d$ to these three equations and
next substituting the expressions of
$d\rho_0,d\zeta_0,d\overline{\zeta}_0$, presented in
\thetag{\ref{d0}}, give:
\begin{equation*}
\footnotesize \aligned \left\{
\begin{array}{l}
d\rho=d{\sf a}\wedge\rho_0+{\sf
a}i\,\zeta_0\wedge\overline{\zeta}_0+{\sf
a}\overline{P}\,\rho_0\wedge\overline{\zeta}_0+{\sf
a}P\,\rho_0\wedge\zeta_0
\\
\\
d\zeta=d{\sf b}\wedge\rho_0+d{\sf c}\wedge\zeta_0+{\sf
b}i\,\zeta_0\wedge\overline{\zeta}_0+{\sf
b}\overline{P}\,\rho_0\wedge\overline{\zeta}_0+{\sf
b}P\,\rho_0\wedge\zeta_0
\\
\\
d\overline{\zeta}=d\overline{{\sf b}}\wedge\rho_0+d\overline{{\sf
c}}\wedge\overline{\zeta}_0+\overline{{\sf
b}}i\,\zeta_0\wedge\overline{\zeta}_0+\overline{{\sf
b}}\overline{P}\,\rho_0\wedge\overline{\zeta}_0+\overline{{\sf
b}}P\,\rho_0\wedge\zeta_0,
\end{array}
\right.
\endaligned
\end{equation*}
or equivalently in matrix notation:
\begin{equation}\footnotesize
\label{dlifted} d\left(%
\begin{array}{c}
\rho \\
\zeta \\
\overline{\zeta} \\
\end{array}%
\right)=
\underbrace{\left(%
\begin{array}{ccc}
d{\sf a} & 0 & 0 \\
d{\sf b} & d{\sf c} & 0 \\
d\overline{{\sf b}} & 0 & d\overline{\sf c} \\
\end{array}%
\right)}_{dg} \wedge
\left(%
\begin{array}{c}
\rho_0 \\
\zeta_0 \\
\overline{\zeta}_0 \\
\end{array}%
\right) +
\left(%
\begin{array}{ccc}
{\sf a}\,P & {\sf a}\,\overline{P} & {\sf a}\,i \\
{\sf b}\,P & {\sf b}\,\overline{P} & {\sf b}\,i \\
\overline{{\sf b}}\,P & \overline{{\sf b}}\,\overline{P} & \overline{{\sf b}}i \\
\end{array}%
\right)
\left(%
\begin{array}{c}
\rho_0\wedge\zeta_0 \\
\rho_0\wedge\overline{\zeta}_0 \\
\zeta_0\wedge\overline{\zeta}_0 \\
\end{array}%
\right).
\end{equation}
On the other hand, multiplying both sides
of~\thetag{ \ref{lifted-coframe}} by the inverse matrix:
\begin{eqnarray*}\footnotesize
g^{-1}=\left(%
\begin{array}{ccc}
\frac{1}{\sf a} & 0 & 0 \\
-\frac{\sf b}{\sf ac} & \frac{1}{\sf c} & 0 \\
-\frac{\overline{\sf b}}{{\sf a}\overline{\sf c}} & 0 & \frac{1}{\overline{\sf c}} \\
\end{array}%
\right)
\end{eqnarray*}
yields the expressions of
$\rho_0,\zeta_0,\overline{\zeta}_0$
in terms of
$\rho,\zeta,\overline{\zeta}$:
\begin{equation}\footnotesize
\label{rho0-rho} \aligned \rho_0&=\frac{1}{\sf a}\rho
\\
\zeta_0&=-\frac{\sf b}{\sf ac}\rho+\frac{1}{\sf c}\zeta
\\
\overline{\zeta}_0&=-\frac{\overline{\sf b}}{{\sf a}\overline{\sf
c}}\rho+\frac{1}{\overline{\sf c}}\overline{\zeta}.
\endaligned
\end{equation}
We may then compute the three exterior products between
these basic $1$-forms:
\begin{equation}
\footnotesize \aligned \label{theta-rep} \left[
\begin{array}{l}
\rho_0\wedge\zeta_0=\frac{1}{{\sf ac}}\rho\wedge\zeta
\\
\\
\rho_0\wedge\overline{\zeta}_0=\frac{1}{{\sf a}\overline{\sf
c}}\rho\wedge\overline{\zeta}
\\
\\
\zeta_0\wedge\overline{\zeta}_0=\frac{\overline{\sf b}}{{\sf
ac}\overline{\sf c}}\rho\wedge\zeta-\frac{\sf b}{{\sf
ac}\overline{\sf c}}\rho\wedge\overline{\zeta}+\frac{1}{{\sf
c}\overline{\sf c}}\zeta\wedge\overline{\zeta}.
\end{array}
\right.
\endaligned
\end{equation}
In addition, one has to replace the first part $dg\wedge
(\rho_0,\zeta_0,\overline\zeta_0)^t$ in \thetag{\ref{dlifted}} by:
\begin{eqnarray}\footnotesize
\label{MC} \underbrace{dg\cdot
g^{-1}}_{\omega_{MC}}\wedge\,\underbrace{
g\,.\,\left(%
\begin{array}{c}
\rho_0 \\
\zeta_0 \\
\overline{\zeta}_0 \\
\end{array}%
\right)}_{(\rho,\zeta,\overline{\zeta})^t},
\end{eqnarray}
and finally we obtain from~\thetag{\ref{dlifted}}, the exterior
differentiations of the lifted 1-forms $\rho,\zeta,\overline\zeta$:
\begin{equation}
\footnotesize\aligned
\label{structure-equation1} d \left(%
\begin{array}{c}
\rho \\
\zeta \\
\overline{\zeta} \\
\end{array}%
\right)&=&\underbrace{
\left(%
\begin{array}{ccc}
\gamma & 0 & 0 \\
\\
\beta & \alpha & 0 \\
\\
\overline{\beta} & 0 & \overline{\alpha} \\
\end{array}%
\right)}_{\omega_{MC}} \wedge
\left(%
\begin{array}{c}
\rho \\
\zeta \\
\overline{\zeta} \\
\end{array}%
\right) +
\left(%
\begin{array}{c}
U_1\,\rho\wedge\zeta
+
\overline U_1\,\rho\wedge\overline{\zeta}
+
U_2\zeta\wedge\overline{\zeta}
\\
\\
V_1\,\rho\wedge\zeta
+
V_2\,\rho\wedge\overline{\zeta}
+
V_3\,\zeta\wedge\overline{\zeta}
\\
\\
\overline V_2\,\rho\wedge\zeta
+
\overline V_1\,\rho\wedge\overline{\zeta}
-
\overline V_3\,\zeta\wedge\overline{\zeta} \\
\end{array}%
\right),
\endaligned
\end{equation}
which incorporate the following {\sl torsion coefficients}:
\begin{eqnarray*}
\begin{array}{ccc}
U_1:=\frac{P\,\overline{\sf c}+\overline{\sf b}i}{{\sf
c}\overline{\sf c}} & U_2:=\frac{{\sf a}i}{{\sf c}\overline{\sf c}}
\\
\ \ V_1:=\frac{P\,{\sf b}\overline{\sf c}+\overline{\sf b}{\sf
b}i}{{\sf a}{\sf c}\overline{\sf c}} & \ \ \ \ \ \ \
V_2:=\frac{\overline{P}\,{\sf b}{\sf c}-{\sf b}^2i}{{\sf a}{\sf
c}\overline{\sf c}} & \ \ \ \ \ \ \ V_3:=\frac{{\sf b}i}{{\sf
c}\overline{\sf c}}
\end{array},
\end{eqnarray*}
and in which the three plain {\sl Maurer-Cartan} 1-forms are:
\[
\aligned \alpha:=\frac{d{\sf c}}{{\sf c}},\ \ \ \ \beta:=\frac{d{\sf
b}}{{\sf a}}-\frac{{\sf b}\,d{\sf c}}{{\sf ac}}, \ \ \ \
\gamma:=\frac{d{\sf a}}{{\sf a}}.
\endaligned
\]
Here the obtained equations are called the {\sl structure equations}
of the problem and moreover the appearing matrix $\omega_{MC}$ is the
so-called {\sl Maurer-Cartan form} of $G$.

\subsection{Absorbtion and normalization}

One of the most essential parts of the Cartan (equivalence) algorithm
is the absorbtion-normalization step, which,
generally speaking, is expressed as follows.

\begin{Observation} (see \cite{5-cubic})
\label{prop-changes} Let $\Theta:=\{\theta^1,\ldots,\theta^n\}$ be a
lifted coframe associated to an equivalence problem having
structure equations{\em :}
\begin{equation*}
\aligned d\theta^i =
\sum_{k=1}^n\bigg(\sum_{s=1}^r\,a_{ks}^i\,\alpha^s +
\sum_{j=1}^{k-1}T^i_{jk}\,\theta^j\bigg)\wedge\theta^k \ \ \ \ \ \ \
\ \ \ \ \ \ {\scriptstyle{(i\,=\,1\,\cdots\,n)}}.
\endaligned
\end{equation*}
Then, one can replace each Maurer-Cartan form $\alpha^s$ and each
torsion coefficient $T^i_{jk}$ with:
\begin{equation}
\label{possible-changes} \boxed{\aligned \alpha^s & \longmapsto
\alpha^s+\sum_{j=1}^n\,z^s_j\,\theta^j \ \ \ \ \ \ \ \ \ \ \ \ \ \ \
\ \ \ \ \ \ \ \ \ \ \ \ \ \ {\scriptstyle{(s\,=\,1\,\cdots\,r)}},
\\
T^i_{jk}&\longmapsto T^i_{jk} + \sum_{s=1}^r\, \big( a_{js}^i\,z_k^s
- a_{ks}^i\,z_j^s \big) \ \ \ \ \ \ \ \ \ \ \
{\scriptstyle{(i\,=\,1\,\cdots\,n\,;\,\,\,
1\,\leqslant\,j\,<\,k\,\leqslant\,n)}},\,
\endaligned}
\end{equation}
for some arbitrary functions $z^\bullet_\bullet$ on the base manifold
$M$. \qed
\end{Observation}

Then one does such a replacement so as to annihilate as many torsion
coefficients as possible, by some appropriate determinations of the
functions $z^\bullet_\bullet$.

Thus, let us perform the following replacements:
\begin{equation}
\label{replacements}
\aligned
\alpha&\mapsto\alpha+p_1\,\rho+q_1\,\zeta+r_1\,\overline\zeta,
\\
\beta&\mapsto\beta+p_2\,\rho+q_2\,\zeta+r_2\,\overline\zeta,
\\
\gamma&\mapsto\gamma+p_3\,\rho+q_3\,\zeta+r_3\,\overline\zeta.
\endaligned
\end{equation}
These substitutions convert the structure equations
\thetag{\ref{structure-equation1}} into the form\,\,---\,\,from now
on and for brevity, we drop presenting the structure equation
$d\overline{\zeta}$ since it is just the conjugation of $d\zeta$:
\[\small
\aligned
d\rho&=\gamma\wedge\rho+(U_1-q_3)\,\rho\wedge\zeta+(\overline
U_1-r_3)\,\rho\wedge\overline\zeta+ U_2\,\zeta\wedge\overline{\zeta},
\\
d\zeta&=\beta\wedge\rho+\alpha\wedge\zeta+(V_1-q_2+p_1)\,\rho\wedge\zeta+(V_2-r_2)\,
\rho\wedge\overline{\zeta}+(V_3-r_1)\,\zeta\wedge\overline{\zeta}.
\endaligned
\]
Visually, one sees that by some appropriate determinations of
$p_i,q_i,r_i$, one can annihilate all the (so modified) torsion
coefficients, except just one, namely $U_2$ in front of
$\zeta\wedge\overline{\zeta}$ at the end of the first line.
Consequently, this torsion coefficient $U_2$ is {\sl essential}, and
the general theory (\cite{ Olver}) shows that $U_2$ (potentially)
provides a {\it normalization} of some group parameter, and here
because $U_2$ is so simple, normalizing it to be $U_2 := i$ provides
the simple group parameter reduction:
\[
\boxed{{\sf a}:={\sf c}\,\overline{\sf c}.}
\]

This then replaces the Maurer-Cartan form $\gamma=\frac{d\sf a}{\sf a}$ by
$\alpha+\overline\alpha$ and transforms the structure equations
\thetag{\ref{structure-equation1}} into the form:
\begin{equation}\aligned
\label{structure-equation} d\rho&=(\alpha+\overline
\alpha)\wedge\rho
+
U_1\,\rho\wedge\zeta+\overline
U_1\,\rho\wedge\overline\zeta+i\,\zeta\wedge\overline\zeta,
\\
d\zeta&=\beta\wedge\rho+\alpha\wedge\zeta+V_1\,\rho\wedge\zeta+V_2\,\rho\wedge\overline\zeta+V_3\,\zeta\wedge\overline\zeta,
\endaligned
\end{equation}
with new torsion coefficients:
\begin{eqnarray*}
\begin{array}{ccc}
U_1:=\frac{P\,\overline{\sf c}+\overline{\sf b}i}{{\sf
c}\overline{\sf c}} &
\\
\ \ V_1:=\frac{P\,{\sf b}\overline{\sf c}+\overline{\sf b}{\sf
b}i}{{\sf c}^2\overline{\sf c}^2} & \ \ \ \ \ \ \
V_2:=\frac{\overline{P}\,{\sf b}{\sf c}-{\sf b}^2i}{{\sf
c}^2\overline{\sf c}^2} & \ \ \ \ \ \ \ V_3:=\frac{{\sf b}i}{{\sf
c}\overline{\sf c}}
\end{array}
\end{eqnarray*}
and with the new Maurer-Cartan 1-forms:
\[
\aligned \alpha:=\frac{d{\sf c}}{{\sf c}} \ \ \ \ \ \
\beta:=\frac{d{\sf b}}{{\sf c}\overline{\sf c}}-\frac{{\sf b}\,d{\sf
c}}{{\sf c}^2\overline{\sf c}}.
\endaligned
\]

Now, let us try again a second absorbtion-normalization procedure.
Doing similar replacements:
\[
\aligned
\alpha&\mapsto\alpha+p_1\,\rho+q_1\,\zeta+r_1\,\overline\zeta,
\\
\beta&\mapsto\beta+p_2\,\rho+q_2\,\zeta+r_2\,\overline\zeta,
\endaligned
\]
one obtains:
\begin{equation}
\aligned\label{absorption} d\rho
&=(\alpha+\overline\alpha)\wedge\rho+(U_1-q_1-\overline
r_1)\rho\wedge\zeta+ (\overline U_1-r_1-\overline
q_1)\rho\wedge\overline{\zeta}+ i\,\zeta\wedge\overline{\zeta},
\\
d\zeta&=\beta\wedge\rho+\alpha\wedge\zeta+(V_1-q_2+p_1)\rho\wedge\zeta+(V_2-
r_2)\rho\wedge\overline{\zeta}+(V_3- r_1)\zeta\wedge\overline{\zeta}.
\endaligned
\end{equation}
Visually, one can annihilate all the (so modified) torsion
coefficients by choosing:
\[\aligned
q_1&:=U_1-\overline V_3, \ \ \ \ \ \ \ \ \ \ r_1:=V_3,
\\
q_2&:=V_1+p_1, \ \ \ \ \ \ \ \ \ \ \ \ r_2:=V_2,
\endaligned
\]
while the two remaining functions:
\[
p_1=:{\sf s}, \ \ \ \ \ \ \ \ \ \ \ \ \ \ p_2=:{\sf r}
\]
can yet be chosen {\it arbitrarily}.

Chosing first these last two functions to be $0$, and
coming back to the explicit expressions of
$U_1$, $V_1$, $V_2$, $V_3$, we see by introducing
the following two {\sl modified Maurer Cartan forms}:
\begin{equation}
\label{alpha0-beta0} \footnotesize\aligned \alpha_0&=\frac{d\sf
c}{\sf c}-\frac{P\overline{\sf c}+2\,i\,\overline{\sf b}}{{\sf
c}\overline{\sf c}}\zeta-\frac{i\sf b}{{\sf c}\overline{\sf
c}}\overline\zeta,
\\
\beta_0&=\frac{d\sf b}{{\sf c}\overline{\sf c} }-\frac{{\sf b}d\sf
c}{{\sf c}^2\overline{\sf c}}-\frac{P{\sf b}\overline{\sf c}+i\,{\sf
b}\overline{\sf b}}{{\sf c}^2\overline{\sf c}^2}\zeta-\frac{\overline
P{\sf bc}-i\,{\sf b}^2}{{\sf c}^2\overline{\sf c}^2}\overline\zeta,
\endaligned
\end{equation}
that the whole torsion is absorbed so that the structure
equations receive the very simple form:
\begin{equation}
\label{struc-eq-after-prolongation-0} \aligned
d\rho&=(\alpha_0+\overline\alpha_0)\wedge\rho+i\,\zeta\wedge\overline\zeta,
\\
d\zeta&=\beta_0\wedge\rho+\alpha\wedge\zeta.
\endaligned
\end{equation}
At this stage, no torsion coefficient can be used anymore
to reduce the structure group.

In fact, one verifies that the
two complex parameters $\sf r$ and ${\sf s}$ and their conjugations
are precisely the {\sl free variables} in the
absorption equations, and consequently, according to the
general procedure, one has to
{\sl prolong} the equivalence problem.

\section{Prolongation of the equivalence problem}
\label{Prolongation}

\subsection{Prolongation procedure}

If one therefore encodes the general remaining
ambiguity in the choice of $\alpha_0$ and $\beta_0$ by
setting:
\begin{equation}
\aligned
\begin{array}{l}\label{modified-MC} \alpha:=\alpha_0+{\sf
s}\rho,
\\
\beta:=\beta_0+{\sf r}\,\rho+{\sf s}\,\zeta,
\end{array}
\endaligned
\end{equation}
one will still have that the absorbed equations look the same
(without lower index `${}_0$'):
\begin{equation}
\boxed{\,
\label{struc-eq-after-prolongation} \aligned
d\rho&=(\alpha+\overline\alpha)\wedge\rho+i\,\zeta\wedge\overline\zeta,
\\
d\zeta&=\beta\wedge\rho+\alpha\wedge\zeta.
\endaligned\,}
\end{equation}

At this moment, one has to launch the prolongation procedure. This
part of Cartan's algorithm relies on the following general result
({\it see} \cite{Olver}, page 395 Proposition 12.13):

\begin{Proposition}
\label{prop.prolongation}
Let $\Theta$ and $\Theta'$ be lifted coframes
of an equivalence problem which admits a non-involutive system of
structure equations and which has a positive degree of indeterminancy.
Let $\Lambda$ and $\Lambda'$ be the modified Maurer-Cartan forms
after the last absorbtion-normalization step. Then, there exists a
diffeomorphism $\Phi:M\longrightarrow M'$ mapping $\Theta$ to
$\Theta'$ for some choice of the group parameters if and only if
there is a diffeomorphism $\Psi:M\times G\longrightarrow M'\times G'$
mapping the coframe $(\Theta,\Lambda)$ to $(\Theta',\Lambda')$ for
some choice of the prolonged group parameters.
\end{Proposition}

This permits us to change our concentration on the original
equivalence problem of the three dimensional hypersurfaces $M^3$
equipped with the lifted coframes $\{\rho,\zeta,\overline\zeta\}$ to
that, along the same lines, of the {\sl prolonged manifolds} $M^{\sf
pr}:=M^3\times G$ with the lifted coframe\,\,---\,\,living on the
product $M^{\sf pr}\times G^{\sf pr}=(M^3\times G)\times G^{\sf
pr}$\,\,---\,\,of the seven 1-forms $\rho,\zeta,\overline{\zeta},
\psi,\varphi,\overline{\psi},\overline{\varphi}$, defined as follows:
\begin{equation}\footnotesize
\label{pro-matrix}
\left(%
\begin{array}{c}
\rho \\
\zeta \\
\overline{\zeta} \\
\alpha \\
\beta \\
\overline{\alpha} \\
\overline{\beta} \\
\end{array}%
\right) =
\underbrace{\left(%
\begin{array}{ccccccc}
1 & 0 & 0 & 0 & 0 & 0 & 0 \\
0 & 1 & 0 & 0 & 0 & 0 & 0 \\
0 & 0 & 1 & 0 & 0 & 0 & 0 \\
{\sf s} & 0 & 0 & 1 & 0 & 0 & 0 \\
{\sf r} & {\sf s} & 0 & 0 & 1 & 0 & 0 \\
\overline{\sf s} & 0 & 0 & 0 & 0 & 1 & 0 \\
\overline{\sf r} & 0 & \overline{\sf s} & 0 & 0 & 0 & 1 \\
\end{array}%
\right)}_{g_{\sf pr}\in G^{\sf pr}}\cdot
\left(%
\begin{array}{c}
\rho \\
\zeta \\
\overline{\zeta} \\
\alpha_0 \\
\beta_0 \\
\overline{\alpha}_0 \\
\overline{\beta}_0 \\
\end{array}%
\right),
\end{equation}
with the new structure group $G^{\sf pr}$, a subgroup of ${\sf
GL}_{3+4}(\mathbb C)={\sf GL}_7(\mathbb C)$ constituted by the {\sl
prolonged group parameters} ${\sf r,s}$ and their conjugates.

\begin{Remark}
\label{Remark}
This prolonged group \thetag{\ref{pro-matrix}} resembles much the
equations~(3.3) on page 7 of the paper~\cite{ GTW}, devoted to the
equivalence problem for second order ordinary differential
equations. In fact, there exists for known reasons ({\em cf.} {\em
e.g.}~\cite{Sparling-Nurowski, Merker}), a certain {\sl transfert
principle} showing that these two seemingly different equivalence
problems will follow fairly the same lines of resolution.  Our main
goal here is to go beyond the so-called\,\,---\,\,usually less
costful\,\,---\,\,{\sl non-parametric} approach and to perform all
computations effectively in terms of the single function $P$, hence in
terms of the graphing function $\varphi(z,\overline z,u)$ of our
hypersurface. In fact, with our choice $\{\mathcal
L,\overline{\mathcal L},\mathcal T\}$ of an initial frame for $TM^3$,
which is explicit in terms of $\varphi$, we deviate from the common
approaches.
\end{Remark}

With the obtained four supplementary $1$-forms
$\alpha$, $\beta$, $\overline{ \alpha}$,
$\overline{ \beta}$, we can now start
the first loop of absorbtion and normalization on the
$7$-dimensional prolonged
space.

\smallskip

Letting a group element
$g_{\sf pr}\in G^{\sf pr}$ be in \thetag{\ref{pro-matrix}}
and abbreviating:
\[
\Omega_0
:=
\big(\rho,\zeta,\overline\zeta,\alpha_0,\beta_0,
\overline\alpha_0,\overline\beta_0\big),
\ \ \ \ \ \ \ \ \ \ \ \ \ \ \ \ \
\ \ \ \ \ \ \ \ \ \ \ \ \ \ \ \ \
\Omega
:=
\big(
\rho,\zeta,\overline\zeta,\alpha,\beta,\overline\alpha,\overline\beta
\big)
\]
the first simple
computation shows that the associated structure equations:
\[
d\Omega
=
\big(dg_{\sf pr}\cdot g_{\sf pr}^{-1}\big)\wedge\Omega+g_{\sf
pr}\cdot d\Omega_0,
\]
read as:
\begin{equation}
\label{d-pr1}
\footnotesize\aligned d\left(%
\begin{array}{c}
\rho \\
\zeta \\
\overline\zeta \\
\alpha \\
\beta \\
\overline\alpha \\
\overline\beta \\
\end{array}%
\right)
=\left(%
\begin{array}{ccccccc}
0 & 0 & 0 & 0 & 0 & 0 & 0 \\
0 & 0 & 0 & 0 & 0 & 0 & 0 \\
0 & 0 & 0 & 0 & 0 & 0 & 0 \\
\delta & 0 & 0 & 0 & 0 & 0 & 0 \\
\gamma & \delta & 0 & 0 & 0 & 0 & 0 \\
\overline\delta & 0 & 0 & 0 & 0 & 0 & 0 \\
\overline\gamma & 0 &\overline\delta & 0 & 0 & 0 & 0 \\
\end{array}%
\right)
\wedge
\left(%
\begin{array}{c}
\rho \\
\zeta \\
\overline\zeta \\
\alpha \\
\beta \\
\overline\alpha \\
\overline\beta \\
\end{array}%
\right)
+
\left(%
\begin{array}{c}
d\rho \\
d\zeta \\
d\overline\zeta \\
{\sf s}\,d\rho+d\alpha_0 \\
{\sf r}\,d\rho+{\sf s}\,d\zeta+d\beta_0 \\
\overline{\sf s}\,d\rho+d\overline\alpha_0 \\
\overline{\sf r}\,d\rho+\overline{\sf s}\,d\overline\zeta+d\overline\beta_0\\
\end{array}%
\right)
\endaligned
\end{equation}
for two new basic Maurer-Cartan 1-forms:
\[
\gamma:=d{\sf r}, \ \ \ \ \ \ \ \ \ \ \ \ \ \ \ \ \ \ \delta:=d\sf s.
\]

To explicitly find the torsion coefficients which should come from the
last four rows of the rightmost $7 \times 1$ matrix, one needs to
express the exterior derivations of $\alpha_0$ and $\beta_0$ in terms
of the lifted 1-forms, and this task is costful, computationally
speaking.  Instead of performing this directly, let us at first employ
a well-know indirect tool ({\it cf.} \cite{Gardner, Olver}) which
temporarily bypasses this computational obstacle and has the virtue of
enabling one to better predict the way the final structure equations
will look like after absorption.

\medskip\noindent{\bf Cartan's (elementary) Lemma.}
\label{Cartan's Lemma}
{\em Let $\{\omega^1,\ldots,\omega^k\}$ be a set of
linearly independent local 1-forms on some
manifold. Then, $k$ arbitrary 1-forms
$\theta^1,\ldots,\theta^k$ satisfy
$\sum_{i=1}^k\theta^i\wedge\omega^i=0$ if and only if they express
$\theta^i=\sum_{j=1}^k A^i_j\, \omega^j$ for some symmetric matrix of local
functions with $A^i_j=A^j_i$.\qed}

\medskip

The truth here is that one intentionally leaves aside
the question of how
these $A_j^i$ could be expressed in terms
of $\theta^1, \dots, \theta^k$, $\omega^1, \dots, \omega^k$.

\smallskip

Now, using the standard
differentiation formula for the exterior product
of two $1$-forms $\lambda$ and $\mu$
(mind the minus sign!):
\[
d\big(\lambda\wedge\mu\big)
=
d\lambda\wedge\mu
-
\lambda
\wedge
d\mu,
\]
the differentiation of the two equations
\thetag{\ref{struc-eq-after-prolongation}}
gives:
\begin{equation}\footnotesize
\label{diff} \left\{
\begin{array}{l}
d^2\rho=0\equiv\big(\underbrace{(d\alpha+2\,i\,\overline{\beta}\wedge\zeta+i\,\beta\wedge\overline{\zeta})+
(d\overline{\alpha}-2\,i\,\beta\wedge\overline{\zeta}-i\,\overline{\beta}
\wedge\zeta)}_{=:\Xi_1}\big)\wedge\rho,
\\
d^2\zeta=0\equiv(\underbrace{d\alpha+2\,i\,\overline{\beta}
\wedge\zeta+i\,\beta\wedge\overline{\zeta}}_{=:\Xi_2})\wedge\zeta+
(\underbrace{d\beta-\beta\wedge\overline{\alpha}}_{=:\Xi_3})\wedge\rho,
\end{array}
\right.
\end{equation}
noticing as a `trick' that the redundant term $2i\, \overline{\beta}
\wedge \zeta$ in $\Xi_1$ helps us to insure the reality relation:
\[
\Xi_1
=
\Xi_2+\overline\Xi_2,
\]
which will be useful for our next:

\begin{Proposition}
\label{Prop-modofoed-MC} The exterior differentials of the
new prolonged lifted $1$-forms $\alpha$ and $\beta$ can be read as:
\begin{eqnarray}
\label{d-modified-MC} \left\{
\begin{array}{l}
d\alpha=\delta^{\sf modified}\wedge\rho+
\\
\ \ \ \ \ \ \ \ \ \ \
+2\,i\,\zeta\wedge\overline{\beta}+i\,\overline{\zeta}\wedge\beta+
W\,\zeta\wedge\overline{\zeta}
\\
\\
d\beta=\gamma^{\sf modified}\wedge\rho+\delta^{\sf
modified}\wedge\zeta+
\\
\ \ \ \ \ \ \ \ \ \ +\beta\wedge\overline{\alpha}
\end{array}
\right.
\end{eqnarray}
for a certain torsion coefficient $W$
which is {\em real}, and for some two
modified Maurer-Cartan $1$-forms
$\delta^{\sf modified}$ and $\gamma^{\sf modified}$.
\end{Proposition}

\proof
Applying Cartan's
Lemma~\ref{Cartan's Lemma}
to~\thetag{\ref{diff}} brings the following expressions of
$\Xi_1,\Xi_2,\Xi_3$ for some three 1-forms $\mathcal A_{ij},\mathcal
B_{ij}$, $\mathcal C$:
\begin{eqnarray*}
\Xi_1&=&-\mathcal C\wedge\rho,
\\
\Xi_2&=&\mathcal A_{11}\wedge\zeta+\mathcal A_{12}\wedge\rho,
\ \ \ \ \ \ \ \ \ \ \ \ \ \ \ \
\Xi_3=\mathcal A_{12}\wedge\zeta+\mathcal A_{22}\wedge\rho.
\end{eqnarray*}
The relation $\Xi_2+\overline\Xi_2-\Xi_1=0$ we
`trickily' insured then reads as:
\[
\mathcal A_{11}\wedge\zeta+\mathcal
B_{11}\wedge\overline{\zeta}+(\mathcal A_{12}+\overline{\mathcal
A}_{12}+\mathcal C)\wedge\rho\equiv 0.
\]
Again, a further application of Cartan's Lemma yields the
(non-explicit) expressions:
\begin{eqnarray*}\footnotesize
\left[
\begin{array}{l}
\mathcal A_{11}=R_{11}\zeta+R_{12}\overline{\zeta}+R_{13}\rho,
\\
\\
\mathcal B_{11}=R_{12}\zeta+R_{22}\overline{\zeta}+R_{23}\rho,
\\
\\
\mathcal A_{12}+\overline{\mathcal A}_{12}+\mathcal
C=R_{13}\zeta+R_{23}\overline{\zeta}+R_{33}\rho,
\end{array}
\right.
\end{eqnarray*}
by means of some complex functions
$R_{ij}, \, i,j=1,2,3$. If we now denote
the two
1-forms $\mathcal A_{12}$ and $\mathcal A_{22}$ by $\delta^{\sf
modified}$ and $\gamma^{\sf modified}$ (respectively), then the
expressions of $\Xi_1,\Xi_2,\Xi_3$ change into:
\begin{equation*}\footnotesize
\aligned
\Xi_1&=\overline{\delta^{\sf modified}}\wedge\rho+\delta^{\sf
modified}\wedge\rho-R_{13}\zeta\wedge\rho-R_{23}\overline{\zeta}\wedge\rho,
\\
\Xi_2&=R_{12}\overline{\zeta}\wedge\zeta+R_{13}\rho\wedge\zeta+\delta^{\sf
modified}\wedge\rho,
\\
\Xi_3&=\delta^{\sf modified}\wedge\zeta+\gamma^{\sf
modified}\wedge\rho.
\endaligned
\end{equation*}
Comparing with the initial expressions of $\Xi_2$, $\Xi_3$,
$\Xi_3$ in \thetag{\ref{diff}} implies that:
\begin{eqnarray}\footnotesize
\label{before-translating} \left\{
\begin{array}{l}
d\alpha=-2\,i\,
\overline{\beta}\wedge\zeta-i\,\beta\wedge\overline{\zeta}
+
\big(\delta^{\sf modified}-R_{13}\zeta\big)
\wedge
\rho-R_{12}\zeta\wedge\overline{\zeta},
\\
\\
d\beta=\beta\wedge\overline{\alpha}+\delta^{\sf
modified}\wedge\zeta+\gamma^{\sf modified}\wedge\rho
\\
\\
d\overline\alpha=2i\,\beta\wedge\overline\zeta+i\,\overline\beta\wedge\zeta+\big(\overline{\delta^{\sf
modified}}-R_{23}\overline\zeta
\big)\wedge\rho+R_{12}\zeta\wedge\overline\zeta.
\end{array}
\right.
\end{eqnarray}
Now granted the equality $\overline{d\alpha}=d\overline\alpha$, one
obtains the following equation, after plain simplifications:
\[
\aligned
-\overline{R_{13}}\,\overline\zeta\wedge\rho+\overline
R_{12}\zeta\wedge\overline\zeta=-R_{23}\overline\zeta\wedge\rho+R_{12}\zeta\wedge\overline\zeta.
\endaligned
\]
Taking account of the linearly independency between
$\overline\zeta\wedge\rho$ and $\zeta\wedge\overline\zeta$, one
immediately concludes that:
\[
R_{23}=\overline R_{13}
\ \ \ \ \ \ \ \ \ \ \ \ \ \ \ \ \
\text{\rm and}
\ \ \ \ \ \ \ \ \ \ \ \ \ \ \ \ \
R_{12}=\overline R_{12}.
\]
In other words, $R_{12}$ is a {\it real}
function and also one can replace $R_{23}$ with $\overline R_{13}$ in
the expression of $d\overline\alpha$. Lastly, the equations
\thetag{\ref{before-translating}} can be transformed as follows after
the substitution $\delta^{\sf modified}-R_{13}\zeta\mapsto
\delta^{\sf modified}$ and putting $W:=-R_{12}$:
\begin{eqnarray}\footnotesize
\label{after-translating} \left\{
\begin{array}{l}
d\alpha=-2\,i\,\overline{\beta}\wedge\zeta-i\,\beta\wedge\overline{\zeta}+(\underbrace{\delta^{\sf
modified}-R_{13}\zeta}_{\mapsto \delta^{\sf
modified}})\wedge\rho\underbrace{-R_{12}}_{+W}\zeta\wedge\overline{\zeta},
\\
d\beta=\beta\wedge\overline{\alpha}+(\underbrace{\delta^{\sf
modified}-R_{13}\zeta}_{\mapsto \delta^{\sf modified}})
\wedge\zeta+\gamma^{\sf modified}\wedge\rho,
\\
d\overline\alpha=2\,i\,{\beta}\wedge\overline\zeta+i\,\overline\beta\wedge{\zeta}+(\underbrace{\overline\delta^{\sf
modified}-\overline R_{13}\overline\zeta}_{\mapsto
\overline\delta^{\sf
modified}})\wedge\rho\underbrace{+R_{12}}_{-W}\zeta\wedge\overline{\zeta}.
\end{array}
\right.
\end{eqnarray}
This completes the proof.
\endproof

The two equations \thetag{\ref{d-modified-MC}} (together with
their unwritten conjugates) and the three equations of
\thetag{\ref{struc-eq-after-prolongation}} constitute the new structure
equations of the problem with $\delta^{\sf modified}$ and
$\gamma^{\sf modified}$ as the {\em modified} Maurer-Cartan forms
after maximal absorbtion of torsion. Thus,
thanks to the above (non-explicit) proposition, one has bypassed
some painful computations, keeping
track of some relevant, somewhat sufficient information,
as Cartan usually did in his papers. Nevertheless, we
will present just at the moment the explicit
expressions of $\delta^{\sf modified}$ and
$\gamma^{\sf modified}$.

\smallskip

Before doing this, let us present the following assertion which
permits one to consider some two fixed expressions of $\delta^{\sf
modified}$ and $\gamma^{\sf modified}$, enjoying
\thetag{\ref{d-modified-MC}}.

\begin{Lemma}
Let $\delta^{\sf modified},\gamma^{\sf modified}$ and
$\delta^{\sf modified}_0,\gamma^{\sf modified}_0$ be two couples of
1-forms satisfying both the same equations~\thetag{\ref{d-modified-MC}}:
\[
\footnotesize
\left[
\aligned
d\alpha
&
=
\delta^{\sf modified}\wedge\rho
+
2i\,\zeta\wedge\overline{\beta}
+
i\,\overline{\zeta}\wedge\beta
+
W\,\zeta\wedge\overline{\zeta},
\\
d\beta
&
=
\gamma^{\sf modified}
\wedge
\rho
+
\delta^{\sf modified}
\wedge
\zeta
+
\beta
\wedge
\overline{\alpha}.
\endaligned\right.
\ \ \ \ \ \ \ \ \ \ \ \ \ \ \ \ \
\left[
\aligned
d\alpha
&
=
\delta_0^{\sf modified}\wedge\rho
+
2i\,\zeta\wedge\overline{\beta}
+
i\,\overline{\zeta}\wedge\beta
+
W\,\zeta\wedge\overline{\zeta},
\\
d\beta
&
=
\gamma_0^{\sf modified}
\wedge
\rho
+
\delta_0^{\sf modified}
\wedge
\zeta
+
\beta
\wedge
\overline{\alpha}.
\endaligned\right.
\]
Then necessarily:
\begin{equation}
\label{p,q} \aligned \delta^{\sf modified}&=\delta^{\sf
modified}_0+{\sf p}\,\rho,
\\
\gamma^{\sf modified}&=\gamma^{\sf modified}_0+{\sf p}\,\zeta+{\sf
q}\,\rho,
\endaligned
\end{equation}
for some arbitrary complex functions $\sf p$ and $\sf q$.
\end{Lemma}

\proof
A plain subtraction yields:
\[
\aligned 0&\equiv(\delta^{\sf modified}-\delta^{\sf
modified}_0)\wedge\rho,
\\
0&\equiv(\gamma^{\sf modified}-\gamma^{\sf
modified}_0)\wedge\rho+(\delta^{\sf modified}-\delta^{\sf
modified}_0)\wedge\zeta.
\endaligned
\]
Now, Cartan's lemma applied to the first equation immediately gives
the first equation of~\thetag{\ref{p,q}}.  Putting then this into the
second equation obtained by subtraction yields, again by means of
Cartan's lemma, the conclusion.
\endproof

Next, a straightforward computation provides a general lemma,
unavoidably required when one wants to perform all
computations explicitly.

\begin{Lemma}
The exterior differential:
\[\footnotesize\aligned
dG=\mathcal L(G)\cdot\zeta_0+\overline{\mathcal
L}(G)\cdot\overline\zeta_0+\mathcal T(G)\cdot \rho_0
\endaligned
\]
of some function $G(z,\overline z,u)$ of class at least $\mathcal{
C}^1$ on the base manifold $M\subset\mathbb C^2$ reexpresses, in
terms of the lifted coframe, as:
\begin{equation}
\label{dP}\footnotesize
dG=\bigg(\frac{1}{\sf c}\mathcal
L(G)\bigg)\cdot\zeta+\bigg(\frac{1}{\overline{\sf
c}}\overline{\mathcal L}(G)\bigg)\cdot\overline\zeta+\bigg(-\frac{\sf
b}{{\sf c}^2\overline{\sf c}}\mathcal L(G)-\frac{\overline{\sf
b}}{{\sf c}\overline{\sf c}^2}\overline{\mathcal L}(G)+\frac{1}{{\sf
c}\overline{\sf c}}\mathcal T(G)\bigg)\cdot\rho.
\qed
\end{equation}
\end{Lemma}

Thus, we may now compare
and inspect the two separate expressions of $d\alpha$ in
\thetag{\ref{d-modified-MC}} and \thetag{\ref{d-pr1}}, namely:
\begin{equation}
\label{2-alpha}
\aligned d\alpha&=\delta^{\sf
modified}\wedge\rho+2\,i\,\zeta\wedge\overline{\beta}+i\,\overline{\zeta}\wedge\beta+
W\,\zeta\wedge\overline{\zeta},
\\
d\alpha&=d\alpha_0+\gamma\wedge\rho+{\sf s}\,d\rho.
\endaligned
\end{equation}
Here, we must compute the
differential $d\alpha_0$  of
$\alpha_0$ given in \thetag{\ref{alpha0-beta0}}:
\[\footnotesize
\aligned d\alpha_0 & = \zero{d\big( {\textstyle{\frac{d{\sf c}}{\sf
c}}} \big)} - \bigg( \frac{1}{\sf c}\, dP - P\,\frac{1}{{\sf c}{\sf
c}}\,d{\sf c} + 2i\,\frac{1}{{\sf c}\overline{\sf c}}\,d\overline{\sf
b} - 2i\, \frac{\overline{\sf b}}{{\sf c}{\sf c}\overline{\sf
c}}\,d{\sf c} - 2i\,\frac{\overline{\sf b}}{{\sf c}\overline{\sf
c}\overline{\sf c}}\, d\overline{\sf c} \bigg) \wedge \zeta -
\\
& \ \ \ \ \ \ \ \ \ \ \ \ \ \ \ - \bigg( \frac{1}{\sf c}\,P +
2i\,\frac{\overline{\sf b}}{{\sf c}\overline{\sf c}} \bigg)\,d\zeta -
\bigg( i\,\frac{1}{{\sf c}\overline{\sf c}}\, d{\sf b} - i\,\frac{\sf
b}{{\sf c}{\sf c}\overline{\sf c}}\, d{\sf c} - i\,\frac{{\sf
b}}{{\sf c}\overline{\sf c}\overline{\sf c}} \bigg) \wedge
\overline{\zeta} - i\,\frac{\sf b}{{\sf c}\overline{\sf
c}}\,d\overline{\zeta}.
\endaligned
\]
Now, thanks to the expressions \thetag{\ref{modified-MC}} and
\thetag{\ref{alpha0-beta0}}, one obtains:
\begin{equation}
\label{db-dc}\footnotesize \aligned d{\sf c}&={\sf
c}\,\alpha_0+\frac{P\overline{\sf c}+2\,i\,\overline{\sf
b}}{\overline{\sf c}}\zeta+\frac{i\,\sf b}{\overline{\sf
c}}\overline\zeta
\\
& \ \ \ \ \ \ \ \ \ \ \ \ ={\sf c}\,\alpha-{\sf
cs}\,\rho+\frac{P\overline{\sf c}+2\,i\,\overline{\sf
b}}{\overline{\sf c}}\zeta+\frac{i\,\sf b}{\overline{\sf
c}}\overline\zeta,
\\
d{\sf b}&={\sf c}\overline{\sf c}\beta_0+{\sf
b}\alpha_0+\frac{2\,P{\sf b}\overline{\sf c}+3\,i\,{\sf
b}\overline{\sf b}}{{\sf c}\overline{\sf c}}\zeta+\frac{\overline
P{\sf b}}{\overline{\sf c}}\overline\zeta
\\
&\ \ \ \ \ \ \ \ \ \ \ \ ={\sf b}\,\alpha+{\sf c}\overline{\sf
c}\,\beta-\big({\sf c}\overline{\sf c}{\sf r}+{\sf bs}\big)\,\rho
+\Big(\frac{2\,P{\sf b}\overline{\sf c}+3\,i\,{\sf b}\overline{\sf
b}}{{\sf c}\overline{\sf c}}-{\sf sc}\overline{\sf
c}\Big)\zeta+\frac{\overline P{\sf b}}{\overline{\sf
c}}\overline\zeta.
\endaligned
\end{equation}
These equations together with \thetag{\ref{dP}} and
\thetag{\ref{struc-eq-after-prolongation}}, enable one to transform
the second expression of $d\alpha$ in \thetag{\ref{2-alpha}} into:
\[
\footnotesize\aligned d\alpha &= \bigg\{\bigg(\frac{\sf b}{{\sf
c}^3\overline{\sf c}}\mathcal L(P)+\frac{\overline{\sf b}}{{\sf
c}^2\overline{\sf c}^2}\overline{\mathcal L}(P)-\frac{1}{{\sf
c}^2\overline{\sf c}}\mathcal T(P)-\frac{P\sf s}{\sf c}+2i\,{\sf
r}-2i\,\frac{{\sf s}\overline{\sf b}}{{\sf c}\overline{\sf
c}}\bigg)\cdot\rho +
\\
&\ \ \ \ \ \ \ \ \ \ \ \ \ \ \ +\bigg(-\frac{1}{{\sf c}\overline{\sf
c}}\overline{\mathcal L}(P)+i\,\frac{P\sf b}{{\sf c}^2\overline{\sf
c}}-2i\,\frac{\overline P\overline{\sf b}}{{\sf c}\overline{\sf
c}^2}-4\,\frac{{\sf b}\overline{\sf b}}{{\sf c}^2\overline{\sf
c}^2}+i\,{\sf s}\bigg)\cdot\overline\zeta-2i\,\overline\beta\bigg\}\,
\wedge \zeta +
\\
&+\bigg\{\bigg(i\,{\sf r}+i\,\frac{{\sf b}\overline{\sf s}}{{\sf
c}\overline{\sf c}}\bigg)\cdot\rho+\bigg(-i\,\frac{P\sf b}{{\sf
c}^2\overline{\sf c}}+2\frac{{\sf b}\overline{\sf b}}{{\sf
c}^2\overline{\sf c}^2}+i\,\frac{P\sf b}{{\sf c}^2\overline{\sf
c}}-\frac{{\sf b}\overline{\sf b}}{{\sf c}^2\overline{\sf
c}^2}+i\,{\sf s}\bigg)\cdot\zeta+i\,\beta \bigg\} \wedge
\overline{\zeta} +
\\
&+\bigg\{\bigg(-\frac{P}{\sf c}-2i\,\frac{\overline{\sf b}}{{\sf
c}\overline{\sf c}}+{\sf s}\bigg)\cdot\beta+\bigg(-i\frac{\sf b}{{\sf
c}\overline{\sf c}}+{\sf
s}\bigg)\cdot\overline\beta+\gamma\bigg\}\wedge\rho.
\endaligned
\]
Chasing then just the coefficient of $\zeta\wedge\overline\zeta$
in this last (long) expression,
which is the function we called $W$, we
therefore obtain the explicit
expression of this single essential torsion coefficient:
\begin{equation}
\label{W}\aligned W=\frac{1}{{\sf c}\overline{\sf
c}}\,\overline{\mathcal L}(P)-2i\,\frac{\sf b}{{\sf c}^2\overline{\sf
c}}P+2i\,\frac{\overline{\sf b}}{{\sf c}\overline{\sf c}^2}\overline
P+6\,\frac{{\sf b}\overline{\sf b}}{{\sf c}^2\overline{\sf
c}^2}+2\,i\,{\sf s}-2\,i\,\overline{\sf s}.
\endaligned
\end{equation}
Thanks to Lemma~\ref{LP}, one easily realizes that $W$ is a real
function as was already mentioned in Proposition \ref{Prop-modofoed-MC}.

Furthermore, collecting together the coefficients of
$\bullet\wedge\rho$ from these two
expressions of $d\alpha$, one also finds
the explicit expression of $\delta^{\sf modified}$:
\begin{equation}
\label{delta0} \small\aligned \delta^{\sf
modified}&=\bigg(\frac{1}{{\sf c}^2\overline{\sf c}}\mathcal
T(P)-\frac{\sf b}{{\sf c}^3\overline{\sf c}}\mathcal
L(P)-\frac{\overline{\sf b}}{{\sf c}^2\overline{\sf
c}^2}\overline{\mathcal L}(P)+\frac{{\sf s}}{\sf
c}P+2i\,\frac{\overline{\sf b}{\sf s}}{{\sf c}\overline{\sf
c}}-2i\,\overline{\sf r}\bigg)\cdot\zeta+\bigg(i\frac{{\sf
b}\overline{\sf s}}{{\sf c}\overline{\sf c}}-i\,{\sf
r}\bigg)\cdot\overline\zeta+
\\
&+{\sf s}\,\alpha-\bigg(\frac{1}{\sf c}P+2i\,\frac{\overline{\sf
b}}{{\sf c}\overline{\sf c}}\bigg)\cdot\beta+{\sf
s}\,\overline\alpha-i\frac{\sf b}{{\sf c}\overline{\sf
c}}\,\overline\beta+
\\
&+d{\sf s}.
\endaligned
\end{equation}

Likewise, let us consider the two separate expressions:
\begin{equation}
\label{dbeta}
\aligned d\beta&=\gamma^{\sf
modified}\wedge\rho+\delta^{\sf
modified}\wedge\zeta+\beta\wedge\overline\alpha,
\\
d\beta&=d\beta_0+\delta\wedge\rho+{\sf r}d\rho+\gamma\wedge\zeta+{\sf
s}\,d\zeta,
\endaligned
\end{equation}
of $d\beta$ in
\thetag{\ref{d-modified-MC}} and \thetag{\ref{d-pr1}}, with
$d\beta_0$ being the differentiation of $\beta_0$ in
\thetag{\ref{alpha0-beta0}} as follows:
\[
\footnotesize \aligned d\beta_0&=\bigg(-\frac{1}{{\sf c}\overline{\sf
c}^2}d\overline{\sf c}\wedge d{\sf b}+\frac{\sf b}{{\sf
c}^2\overline{\sf c}^2}d\overline{\sf c}\wedge d{\sf
c}\bigg)-\bigg(\frac{P\sf b}{{\sf c}^2\overline{\sf c}}+i\,\frac{{\sf
b}\overline{\sf b}}{{\sf c}^2\overline{\sf c}^2}\bigg)d\zeta+
\\
&+\bigg(-\frac{\sf b}{{\sf c}^2\overline{\sf c}}dP-\frac{P}{{\sf
c}^2\overline{\sf c}}d{\sf b}+\frac{P\sf b}{{\sf c}^2\overline{\sf
c}^2}d\overline{\sf c}+2\,\frac{P\sf b}{{\sf c}^3\overline{\sf
c}}d{\sf c}-i\frac{\overline{\sf b}}{{\sf c}^2\overline{\sf
c}^2}d{\sf b}-i\,\frac{\sf b}{{\sf c}^2\overline{\sf
c}^2}d\overline{\sf b}+2i\,\frac{{\sf b}\overline{\sf b}}{{\sf
c}^3\overline{\sf c}^2}d{\sf c}+2i\,\frac{{\sf b}\overline{\sf
b}}{{\sf c}^2\overline{\sf c}^3}d\overline{\sf c}\bigg)\wedge\zeta-
\\
&+\bigg(-\frac{\overline P\sf b}{{\sf c}\overline{\sf
c}^2}+i\,\frac{{\sf b}^2}{{\sf c}^2\overline{\sf
c}^2}\bigg)d\overline\zeta+
\\
&+\bigg(-\frac{\sf b}{{\sf c}\overline{\sf c}^2}d\overline
P-\frac{\overline P}{{\sf c}\overline{\sf c}^2}d{\sf
b}+\frac{\overline P{\sf b}}{{\sf c}^2\overline{\sf c}^2}d{\sf
c}+2\,\frac{\overline P{\sf b}}{{\sf c}\overline{\sf
c}^3}d\overline{\sf c}+2i\,\frac{\sf b}{{\sf c}^2\overline{\sf
c}^2}d{\sf b}-2i\,\frac{{\sf b}^2}{{\sf c}^3\overline{\sf c}^2}d{\sf
c}-2i\,\frac{{\sf b}^2}{{\sf c}^2\overline{\sf c}^3}d\overline{\sf
c}\bigg)\wedge\overline\zeta.
\endaligned
\]
Performing lines of (rather lengthy) computations similar
to those we
already did, we can extract the coefficients of $\bullet\wedge\rho$
from the two equal expressions of $d\beta$ in \thetag{\ref{dbeta}} and
we find:
\begin{equation}
\label{gamma0} \footnotesize\aligned \gamma^{\sf
modified}&=\bigg(\frac{\sf b}{{\sf c}^3\overline{\sf c}^2}\mathcal
T(P)-\frac{{\sf b}^2}{{\sf c}^4\overline{\sf c}^2}\mathcal
L(P)-\frac{{\sf b}\overline{\sf b}}{{\sf c}^3\overline{\sf
c}^3}\overline{\mathcal L}(P)+\frac{\sf bs}{{\sf c}^2\overline{\sf
c}}P-\frac{\sf r}{\sf c}P+i\,\frac{{\sf b}\overline{\sf b}\sf s}{{\sf
c}^2\overline{\sf c}^2}-2i\,\frac{\overline{\sf b}\sf r}{{\sf
c}\overline{\sf c}}-i\,\frac{{\sf b}\overline{\sf r}}{{\sf
c}\overline{\sf c}}+{\sf s}\overline{\sf s}\bigg)\cdot\zeta+
\\
&+\bigg(\frac{\sf b}{{\sf c}^2\overline{\sf c}^3}\mathcal T(\overline
P)-\frac{{\sf b}^2}{{\sf c}^3\overline{\sf c}^3}\mathcal L(\overline
P)-\frac{{\sf b}\overline{\sf b}}{{\sf c}^2\overline{\sf
c}^4}\overline{\mathcal L}(\overline P)+\frac{ {\sf b}\overline{\sf
s}}{{\sf c}\overline{\sf c}^2}\overline P-i\,\frac{{\sf
b}^2\overline{\sf s}}{{\sf c}^2\overline{\sf
c}^2}\bigg)\cdot\overline\zeta+
\\
&+{\sf r}\,\alpha-\bigg(\frac{\sf b}{{\sf
c}^2\overline{\sf c}}P+i\,\frac{{\sf b}\overline{\sf b}}{{\sf
c}^2\overline{\sf c}^2}-{\sf s}+\overline{\sf
s}\bigg)\cdot\beta+2{\sf r}\,\overline\alpha+\bigg(-\frac{\sf b}{{\sf
c}\overline{\sf c}^2}\overline P+i\,\frac{{\sf b}^2}{{\sf
c}^2\overline{\sf c}^2}\bigg)\cdot\overline\beta-
\\
&+d{\sf r}.
\endaligned
\end{equation}

From now on and for the sake of simplicity and compatibility among
the notations, let us drop the word "$\sf modified$" from
$\delta^{\sf modified}$ and $\gamma^{\sf modified}$ and denote them
simply by $\delta$ and $\gamma$. Summarizing the results, now the
structure equations \thetag{\ref{d-pr1}} is transformed into:
\begin{eqnarray}
\label{modified-structure-equation}
\begin{array}{ll}
d\rho=\alpha\wedge\rho+\overline\alpha\wedge\rho+i\,\zeta\wedge\overline\zeta,
&
\\
d\zeta=\beta\wedge\rho+\alpha\wedge\zeta, &
\\
d\overline\zeta=\overline\beta\wedge\rho+\overline\alpha\wedge\overline\zeta,
&
\\
d\alpha=\delta\wedge\rho+2\,i\,\zeta\wedge\overline{\beta}+i\,\overline{\zeta}\wedge\beta+

W\,\zeta\wedge\overline{\zeta}, &
\\
d\beta=\gamma\wedge\rho+\delta\wedge\zeta+\beta\wedge\overline{\alpha},
&
\\
d\overline\alpha=\overline\delta\wedge\rho-2\,i\,\overline\zeta\wedge{\beta}-i\,{\zeta}\wedge\overline\beta-

\overline W\,\zeta\wedge\overline{\zeta}, &
\\
d\overline\beta=\overline\gamma\wedge\rho+\overline\delta\wedge\overline\zeta+\overline\beta\wedge{\alpha},
&
\end{array}
\end{eqnarray}
with the already modified Maurer-Cartan forms $\delta$ and $\gamma$
given by \thetag{\ref{delta0}} and \thetag{\ref{gamma0}},
and with some relevant real torsion coefficient $W$ given by
\thetag{\ref{W}}.

\subsection{Absorbtion-normalization}
After having re-shaped so the structure equations, one has to apply
again the absorbtion-normalization procedure by considering the
substitutions:
\[
\aligned \delta&\mapsto
\delta+p_1\,\rho+q_1\,\zeta+r_1\,\overline\zeta+s_1\,\alpha+t_1\,\overline\alpha+u_1\,\beta+v_1\,\overline\beta,
\\
\gamma&\mapsto\gamma+p_2\,\rho+q_2\,\zeta+r_2\,\overline\zeta+s_2\,\alpha+t_2\,\overline\alpha+u_2\,\beta+v_2\,\overline\beta.
\endaligned
\]
One easily verifies by elementary linear algebra computations that
here the single torsion coefficient $W$ is, as guessed, indeed
normalizable.

Normalizing then this coefficient to zero determines $\overline{\sf
s}$ as:
\begin{equation}
\label{s}
\boxed{\,
\overline{\sf s}
=
{\sf s}
-
\frac{i}{2}\frac{1}{{\sf c}\overline{\sf c}}\overline{\mathcal
L}(P)
-
\frac{\sf b}{{\sf c}^2\overline{\sf c}}P
+
\frac{\overline{\sf
b}}{{\sf c}\overline{\sf c}^2}\overline P
-
3i\,\frac{{\sf
b}\overline{\sf b}}{{\sf c}^2\overline{\sf c}^2}.\,}
\end{equation}

Consequently, one has to differentiate this equation:
\[
\footnotesize\aligned d\overline{\sf s}&=d{\sf s}-\bigg\{
3i\,\frac{\sf b}{{\sf c}^2\overline{\sf c}^2}d{\overline{\sf
b}}+3i\,\frac{\overline{\sf b}}{{\sf c}^2\overline{\sf c}^2}d{\sf
b}-6i\,\frac{{\sf b}\overline{\sf b}}{{\sf c}^3\overline{\sf
c}^2}d{\sf c}-6i\,\frac{{\sf b}\overline{\sf b}}{{\sf
c}^2\overline{\sf c}^3}d\overline{\sf c}+\frac{P}{{\sf
c}^2\overline{\sf c}}d{\sf b}+\frac{\sf b}{{\sf c}^2\overline{\sf
c}}dP-2\,\frac{P\sf b}{{\sf c}^3\overline{\sf c}}d{\sf c}-\frac{P\sf
b}{{\sf c}^2\overline{\sf c}^2}d\overline{\sf c}-
\\
&\ \ \ \ \ \ \ \ -\frac{\overline P}{{\sf c}\overline{\sf
c}^2}d\overline{\sf b}-\frac{\overline{\sf b}}{{\sf c}\overline{\sf
c}^2}d\overline P+\frac{\overline P\,\overline{\sf b}}{{\sf
c}^2\overline{\sf c}^2}d{\sf c}+2\frac{\overline P\overline {\sf
b}}{{\sf c}\overline{\sf c}^3}d\overline{\sf c}-\frac{i}{2{\sf
c}^2\overline{\sf c}}\overline{\mathcal L}(P)d{\sf c}-\frac{i}{2{\sf
c}\overline{\sf c}^2}\overline{\mathcal L}(P)d\overline{\sf c}+
+\frac{i}{2{\sf c}\overline{\sf c}}d\overline{\mathcal L}(P) \bigg\},
\endaligned
\]
in which similarly to \thetag{\ref{dP}}, one has:
\begin{equation}
\footnotesize\aligned
\label{dLP} d(\overline{\mathcal L}(P))&=\bigg(\frac{1}{\sf c}\mathcal
L(\overline{\mathcal
L}(P))\bigg)\cdot\zeta+\bigg(\frac{1}{\overline{\sf
c}}\overline{\mathcal L}(\overline{\mathcal
L}(P))\bigg)\cdot\overline\zeta+\bigg(-\frac{\sf b}{{\sf
c}^2\overline{\sf c}}\mathcal L(\overline{\mathcal
L}(P))-\frac{\overline{\sf b}}{{\sf c}\overline{\sf
c}^2}\overline{\mathcal L}(\overline{\mathcal L}(P))+\frac{1}{{\sf
c}\overline{\sf c}}\mathcal T(\overline{\mathcal
L}(P))\bigg)\cdot\rho.
\endaligned
\end{equation}
Then, putting the expressions \thetag{\ref{db-dc}} of $d{\sf b},d{\sf
c}$ into the above equation expression of $d\overline{\sf s}$ changes
it into the following form after simplification:
\[
\footnotesize\aligned d\overline{\sf s}&= d{\sf
s}+
\bigg(-\frac{\overline P\,\overline{\sf
r}}{\overline{\sf c}}+\frac{P\,\sf r}{\sf c}-9\frac{{\sf
b}^2\overline{\sf b}^2}{{\sf c}^4\overline{\sf c}^4}+
\frac{\mathcal
L(P){\sf b}^2}{{\sf c}^4\overline{\sf c}^2}+\frac{PP{\sf b}^2}{{\sf
c}^4\overline{\sf c}^2}+\frac{\overline PP\overline{\sf b}^2}{{\sf
c}^2\overline{\sf c}^4}-\frac{\overline{\sf b}^2\overline{\mathcal
L}(\overline P)}{{\sf c}^2\overline{\sf c}^4}-
\frac{1}{4}\frac{\mathcal L(\overline P)\mathcal L(\overline P)}{{\sf
c}^2\overline{\sf c}^2}+i\frac{\overline{\mathcal L}(\mathcal
L(P))\sf b}{{\sf c}^3\overline{\sf c}^2}+
\\
&+ i\,\frac{\mathcal L(\overline{\mathcal
L}(\overline P))\overline{\sf b}}{{\sf c}^2\overline{\sf
c}^3}-3\frac{\mathcal L(\overline P){\sf b}\overline{\sf b}}{{\sf
c}^3\overline{\sf c}^3}-2\frac{P\sf bs}{{\sf c}^2\overline{\sf c}}-
2\,\frac{P\overline P{\sf b}\overline{\sf b}}{{\sf c}^3\overline{\sf
c}^3}+2\,\frac{\overline P\overline{\sf b}\sf s}{{\sf c}\overline{\sf
c}^2}+i\,\frac{P\mathcal L(\overline P)\sf b}{{\sf c}^3\overline{\sf
c}^2}- i\,\frac{\mathcal L(\overline P)\sf s}{{\sf c}\overline{\sf
c}}-i\,\frac{\overline P\mathcal L(\overline P)\overline{\sf b}}{{\sf
c}^2\overline{\sf c}^3}-
\\
&-\frac{i}{2}\frac{\mathcal L(\mathcal L(\overline P))\sf b}{{\sf
c}^3\overline{\sf c}^2}+ 3i\,\frac{\overline{\sf b}\sf r}{{\sf
c}\overline{\sf c}}-\frac{i}{2}\frac{\mathcal T(\overline{\mathcal
L}(P))}{{\sf c}^2\overline{\sf c}^2}-6i\,\frac{\overline P{\sf
b}\overline{\sf b}^2}{{\sf c}^3\overline{\sf c}^4}+ 3i\,\frac{{\sf
b}\overline{\sf r}}{{\sf c}\overline{\sf c}}-6\,\frac{{\sf
b}\overline{\sf b}\sf s}{{\sf c}^2\overline{\sf c}^2}+6i\,\frac{P{\sf
b}^2\overline{\sf b}}{{\sf c}^4\overline{\sf c}^3}-
\frac{i}{2}\,\frac{\overline{\mathcal L}(\mathcal L(\overline
P))}{{\sf c}^2\overline{\sf c}^3}\bigg)\cdot\rho+
\\
&+ \bigg(\frac{P\sf s}{\sf c}-\frac{i}{2}\frac{\mathcal L(\mathcal
L(\overline P))}{{\sf c}^2\overline{\sf c}}-\frac{\mathcal L(P)\sf
b}{{\sf c}^3\overline{\sf c}}+3i\,\frac{\overline{\sf b}\sf s}{{\sf
c}\overline{\sf c}}+\frac{i}{2}\frac{P\mathcal L(\overline P)}{{\sf
c}^2\overline{\sf c}}+ 3\,\frac{{\sf b}\overline{\sf b}^2}{{\sf
c}^3\overline{\sf c}^3}+\frac{1}{2}\frac{\mathcal L(\overline
P)\overline{\sf b}}{{\sf c}^2\overline{\sf c}^2}-3i\,\frac{P{\sf
b}\overline{\sf b}}{{\sf c}^3\overline{\sf c}^2}\bigg)\cdot\zeta+
\\
&+\bigg(-\frac{\overline P^2\,\overline{\sf b}}{{\sf c}\overline{\sf
c}^3}+6\,\frac{{\sf b}^2\overline{\sf b}}{{\sf c}^3\overline{\sf
c}^3}+ i\,\frac{\mathcal L(\overline P)\,\overline P}{{\sf
c}\overline{\sf c}^2}+\frac{\overline{\mathcal L}(\overline P)}{{\sf
c}\overline{\sf c}^3}+3i\,\frac{\sf bs}{{\sf c}\overline{\sf
c}}-3i\,\frac{P{\sf b}^2}{{\sf c}^3\overline{\sf c}^2}-
\frac{\overline P\sf s}{\overline{\sf c}}+3i\,\frac{\overline P{\sf
b}\overline{\sf b}}{{\sf c}^2\overline{\sf c}^3}+\frac{\mathcal
L(\overline P)\sf b}{{\sf c}^2\overline{\sf c}^2}+\frac{P\overline
P\sf b}{{\sf c}^2\overline{\sf
c}^2}-\frac{i}{2}\frac{\overline{\mathcal L}(\mathcal L(\overline
P))}{{\sf c}\overline{\sf c}^2}\bigg)\cdot\overline\zeta+
\\
&+\bigg(3i\,\frac{{\sf b}\overline{\sf b}}{{\sf c}^2\overline{\sf
c}^2}+\frac{i}{2}\frac{\mathcal L(\overline P)}{{\sf c}\overline{\sf
c}}-\frac{\overline P\,\overline{\sf b}}{{\sf c}\overline{\sf
c}^2}+\frac{P\sf b}{{\sf c}^2\overline{\sf c}}\bigg)\cdot\alpha+
\bigg(-\frac{P}{\sf c}-3i\,\frac{\overline{\sf b}}{{\sf
c}\overline{\sf c}}\bigg)\cdot\beta+\bigg(\frac{3i\,{\sf
b}\overline{\sf b}}{{\sf c}^2\overline{\sf
c}^2}+\frac{i}{2}\frac{\mathcal L(\overline P)}{{\sf c}\overline{\sf
c}}- \frac{\overline P\,\overline{\sf b}}{{\sf c}\overline{\sf
c}^2}+\frac{P\sf b}{{\sf c}^2\overline{\sf
c}}\bigg)\cdot\overline\alpha+\bigg(\frac{\overline P}{\overline{\sf
c}}-3i\,\frac{\sf b}{{\sf c}\overline{\sf
c}}\bigg)\cdot\overline\beta.
\endaligned
\]

Next, by a careful glance on the expression of $\delta$ and its
conjugation ({\it see} \thetag{\ref{delta0}}), we realize that having
$d\overline {\sf s}$ in terms of $d\sf s$ and the lifted 1-forms
$\rho,\zeta,\overline\zeta,\alpha,\beta,\overline\alpha,\overline\beta$
enables us to express $\overline\delta$ in terms of $\delta$ and the
lifted coframe ({\it cf.} \thetag{\ref{delta0}}). More precisely, our
computations show that we have\,\,---\,\,the coefficients of
$\alpha,\beta,\overline\alpha,\overline\beta$ vanish identically
after simplification:
\begin{equation}
\label{deltabar} \aligned
\overline\delta&:=\delta+i\,W_1\,\rho+W_2\,\zeta-\overline
W_2\,\overline\zeta,
\endaligned
\end{equation}
with the coefficients:
\[
\footnotesize\aligned W_1&:=-\frac{1}{2}\frac{\mathcal
T(\overline{\mathcal L}(P))}{{\sf c}^2\overline{\sf
c}^2}+\frac{\mathcal L(\overline{\mathcal L}(\overline
P))\overline{\sf b}}{{\sf c}^2\overline{\sf
c}^3}-\frac{1}{2}\,\frac{\mathcal L(\mathcal (\overline P))\sf
b}{{\sf c}^3\overline{\sf c}^2}-
\frac{1}{2}\,\frac{\overline{\mathcal L}(\mathcal L(\overline
P))\overline{\sf b}}{{\sf c}^2\overline{\sf
c}^3}+\frac{\overline{\mathcal L}(\mathcal L(P))\sf b}{{\sf
c}^3\overline{\sf c}^2} - i\,\frac{\mathcal L(P){\sf b}^2}{{\sf
c}^4\overline{\sf c}^2}+i\,\frac{\overline{\mathcal L}(\overline
P)\overline{\sf b}^2}{{\sf c}^2\overline{\sf c}^4}+
\\
&+\bigg(-\frac{1}{2}\frac{\mathcal L(\overline P)}{{\sf
c}\overline{\sf c}}+i\,\frac{P\sf b}{{\sf c}^2\overline{\sf
c}}-3\,\frac{{{\sf b}\overline{\sf b}}}{{\sf c}^2\overline{\sf
c}^2}-i\,\frac{\overline P\,\overline{\sf b}}{{\sf c}\overline{\sf
c}^2}\bigg)\,{\sf s}+\bigg(3\,\frac{\overline{\sf b}}{{\sf
c}\overline{\sf c}}-i\,\frac{P}{\sf c}\bigg)\,{\sf
r}+\bigg(-\frac{1}{2}\,\frac{\mathcal L(\overline P)}{{\sf
c}\overline{\sf c}}+i\,\frac{P\sf b}{{\sf c}^2\overline{\sf
c}}-3\,\frac{{{\sf b}\overline{\sf b}}}{{\sf c}^2\overline{\sf
c}^2}-i\,\frac{\overline P\,\overline{\sf b}}{{\sf c}\overline{\sf
c}^2}\bigg)\,\overline{\sf s}+\bigg(3\,\frac{\sf b}{{\sf
c}\overline{\sf c}}+i\,\frac{\overline P}{\overline{\sf
c}}\bigg)\,\overline{\sf r} ,
\\
W_2&:=i\,\frac{\overline{\mathcal L}(\mathcal L(P))}{{\sf
c}^2\overline{\sf c}}-\frac{3}{2}i\,\frac{\mathcal L(\mathcal
L(\overline P))}{{\sf c}^2\overline{\sf
c}}+\frac{3}{2}\,\frac{\mathcal L(\overline P)\overline{\sf b}}{{\sf
c}^2\overline{\sf c}^2}+\frac{i}{2}\,\frac{P\mathcal L(\overline
P)}{{\sf c}^2\overline{\sf c}}-3i\,\frac{P {\sf b}\overline{\sf
b}}{{\sf c}^3\overline{\sf c}^2}+ 3\,\frac{{\sf b}\overline{\sf
b}^2}{{\sf c}^3\overline{\sf c}^3}+3i\,\overline{\sf r}.
\endaligned
\]
(We notice {\em passim} that the first torsion coefficient $W_1$ is
real.)

Further, after determining $\overline{\sf s}$ in
\thetag{\ref{s}}, the expressions of $\overline\alpha$ and
$\overline\beta$  change
and are not anymore the conjugates of $\alpha$ and $\beta$. Hence, we
replace the notations $\overline\alpha$ and $\overline\beta$ by
$\widetilde{\alpha}$ and $\widetilde{{\beta}}$, respectively. Putting
this new expression of $\overline\delta$ into the last structure
equation \thetag{\ref{modified-structure-equation}} changes it into
the form:
\begin{eqnarray}
\label{modified-structure-equation-2}
\begin{array}{ll}
d\rho=\alpha\wedge\rho+\widetilde{\alpha}\wedge\rho+i\,\zeta\wedge\overline\zeta,
&
\\
d\zeta=\beta\wedge\rho+\alpha\wedge\zeta, &
\\
d\overline\zeta=\widetilde{\beta}\wedge\rho+\widetilde{\alpha}\wedge\overline\zeta,
&
\\
d\alpha=\delta\wedge\rho+2\,i\,\zeta\wedge\overline{\beta}+i\,\overline{\zeta}\wedge\beta,
&
\\
d\beta=\gamma\wedge\rho+\delta\wedge\zeta+\beta\wedge\overline{\alpha},
 &
 \\
 d\widetilde{\alpha}=\delta\wedge\rho-2\,i\,\overline\zeta\wedge{\beta}-i\,{\zeta}\wedge\widetilde{\beta}+
W_2\,\zeta\wedge\rho-\overline W_2\overline\zeta\wedge\rho, &
 \\
d\widetilde{\beta}=\overline\gamma\wedge\rho+\delta\wedge\overline\zeta+\widetilde{\beta}\wedge{\alpha}+
i\,W_1\,\rho\wedge\overline\zeta+W_2\,\zeta\wedge \overline\zeta.

\end{array}
\end{eqnarray}

\subsection{Absorbtion-normalization of the latest structure equation}
To determine essential torsion coefficients, similarly as before,
we make substitutions of the kind:
\[
\aligned \delta&\mapsto
\delta+p_1\,\rho+q_1\,\zeta+r_1\,\overline\zeta+s_1\,\alpha+t_1\,\overline\alpha+u_1\,\beta+v_1\,\overline\beta,
\\
\gamma&\mapsto\gamma+p_2\,\rho+q_2\,\zeta+r_2\,\overline\zeta+s_2\,\alpha+t_2\,\overline\alpha+u_2\,\beta+v_2\,\overline\beta.
\endaligned
\]
This converts the structure equations into the form:
\[
\footnotesize\aligned
d\alpha&=\delta\wedge\rho+q_1\,\zeta\wedge\rho+r_1\,\overline\zeta\wedge\rho+s_1\,\alpha\wedge\rho+t_1\,\widetilde{\alpha}\wedge\rho+u_1\,\beta\wedge\rho+
v_1\,\widetilde{\beta}\wedge\rho+2i\,\zeta\wedge\widetilde{\beta}+i\,\overline\zeta\wedge\beta,
\\
d\beta&=\gamma\wedge\rho+\delta\wedge\zeta+(q_2-p_1)\zeta\wedge\rho+r_2\,\overline\zeta\wedge\rho+s_2\,\alpha\wedge\rho+t_2\,\widetilde{\alpha}\wedge\rho+
u_2\,\beta\wedge\rho+v_2\,\widetilde{\beta}\wedge\rho+r_1\,\overline\zeta\wedge\zeta+
\\
&\ \ \
+s_1\,\alpha\wedge\zeta+t_1\,\widetilde{\alpha}\wedge\zeta+u_1\,\beta\wedge\zeta+v_1\,\widetilde{\beta}\wedge\zeta+\beta\wedge\widetilde{\alpha},
\\
d\widetilde{\alpha}&=\delta\wedge\rho+(q_1+W_2)\zeta\wedge\rho+(r_1-\overline
W_2)\overline\zeta\wedge\rho+s_1\,\alpha\wedge\rho+t_1\,\widetilde{\alpha}\wedge\rho+u_1\,\beta\wedge\rho+v_1\,\widetilde{\beta}\wedge\rho
-2i\,\overline\zeta\wedge\beta-i\,\zeta\wedge\widetilde{\beta},
\\
d\widetilde{\beta}&=\overline\gamma\wedge\rho+\delta\wedge\overline\zeta+(\overline
q_2-p_1-i\,W_1)\overline\zeta\wedge\rho+\overline
r_2\,\zeta\wedge\rho+\overline
s_2\,\widetilde{\alpha}\wedge\rho+\overline
t_2\,\alpha\wedge\rho+\overline
u_2\,\widetilde{\beta}\wedge\rho+\overline v_2\,\beta\wedge\rho+
\\
&\ \ \
+(q_1+W_2)\zeta\wedge\overline\zeta+s_1\,\alpha\wedge\overline\zeta+t_1\,\widetilde{\alpha}\wedge\overline\zeta+u_1\,\beta\wedge\overline\zeta+
v_1\,\widetilde{\beta}\wedge\overline\zeta.
\endaligned
\]
In order to annihilate as much as possible the appearing (modified)
torsion coefficients, we have to solve the following system of homogeneous
equations:
\[
\aligned  0&=q_1=r_1=s_1=t_1=u_1=v_1, \ \ \ 0=
r_2=s_2=t_2=u_2=v_2,
\\
& 0= q_2-p_1, \ \ \  0= q_1+W_2, \ \ \ 0= r_1-\overline W_2, \ \ \ 0=
\overline q_2-p_1-i\,W_1.
\endaligned
\]
One readily realizes that besides the following determinations:
\[
\aligned q_1&=0, \ \ \ r_i=s_i=t_i=u_i=v_i=0, \ i=1,2,
\\
& \ \ \ \ \ q_2=p_1, \ \ \ \ \ {\rm Im}(p_1)=-\frac{1}{2}W_1,
\endaligned
\]
the homogeneous system will be satisfied if and only if we also have:
\begin{eqnarray*}
\label{normalizable-3} 0\equiv W_2.
\end{eqnarray*}
In other words, $W_2$ is the only normalizable expression of this step. A
careful glance at the expression of this function shows that it will
be normalized to zero as soon as we put:
\begin{equation}
\label{r} \footnotesize\aligned\boxed{{\sf
r}:=-\frac{1}{3}\,\frac{{\mathcal L}(\overline{\mathcal L}(\overline
P))}{{\sf c}\overline{\sf c}^2}+\frac{1}{2}\,\frac{\overline{\mathcal
L}(\overline{\mathcal L}(P))}{{\sf c}\overline{\sf
c}^2}-\frac{i}{2}\,\frac{\overline{\mathcal L}(P){\sf b}}{{\sf
c}^2\overline{\sf c}^2}-\frac{1}{6}\,\frac{\overline
P\,\overline{\mathcal L}(P)}{{\sf c}\overline{\sf
c}^2}+\frac{\overline P {\sf b}\overline{\sf b}}{{\sf
c}^2\overline{\sf c}^3}-i\,\frac{{\sf b}^2\overline{\sf b}}{{\sf
c}^3\overline{\sf c}^3}.}
\endaligned
\end{equation}

With this expression of $\sf r$ which
reduces the group dimension, the only remaining (inessential) torsion
coefficient $W_1$ takes the form:
\begin{equation}\label{W1} \footnotesize\aligned W_1&=-\frac{1}{2}\frac{\mathcal
T(\overline{\mathcal L}(P))}{{\sf c}^2\overline{\sf
c}^2}+\frac{\overline{\mathcal L}(\mathcal L(P))\overline{\sf
b}}{{\sf c}^2\overline{\sf c}^3}- \frac{1}{2}\,\frac{\mathcal
L(\mathcal L(\overline P))\sf b}{{\sf c}^3\overline{\sf
c}^2}+\frac{i}{3}\frac{P\mathcal L(\overline{\mathcal L}(\overline
P))}{{\sf c}^2\overline{\sf c}^2}-\frac{i}{3}\frac{\overline
P\overline{\mathcal L}(\mathcal L(P))}{{\sf c}^2\overline{\sf c}^2}+
\frac{i}{2}\,\frac{\overline P\mathcal L(\overline{\mathcal
L}(P))}{{\sf c}^2\overline{\sf c}^2}-
\\
&- \frac{i}{2}\,\frac{P\overline{\mathcal L}(\mathcal L(\overline
P))}{{\sf c}^2\overline{\sf c}^2}+\frac{3}{2}\frac{\mathcal
L(\overline{\mathcal L}(P))\sf b}{{\sf c}^3\overline{\sf
c}^2}+3i\,\frac{\overline{\mathcal L}(P){{\sf b}\overline{\sf
b}}}{{\sf c}^3\overline{\sf c}^3}+ \frac{i}{6}\frac{P\overline
P\mathcal L(\overline P)}{{\sf c}^2\overline{\sf
c}^2}-\frac{i}{6}\frac{P\overline P\overline{\mathcal L}(P)}{{\sf
c}^2\overline{\sf c}^2}+\frac{i}{4}\frac{\mathcal L(\overline
P)\overline{\mathcal L}(P)}{{\sf c}^2\overline{\sf
c}^2}+i\,\frac{\overline{\mathcal L}(\overline P)\overline{\sf
b}^2}{{\sf c}^2\overline{\sf c}^4}-
\\
&-i\,\frac{\mathcal L(P){\sf b}^2}{{\sf c}^4\overline{\sf
c}^2}-\frac{\overline P\,\overline{\mathcal L}(P)\overline{\sf
b}}{{\sf c}^2\overline{\sf c}^3}-\frac{\overline P\mathcal
L(\overline P)\overline{\sf b}}{{\sf c}^2\overline{\sf
c}^3}+2i\,\frac{P\overline P{\sf b}\overline{\sf b}}{{\sf
c}^3\overline{\sf c}^3}-i\,\frac{\overline P^2\overline{\sf
b}^2}{{\sf c}^2\overline{\sf c}^4}-4\,\frac{\overline P{{\sf
b}\overline{\sf b}^2}}{{\sf c}^3\overline{\sf c}^4}- i\,\frac{P^2{\sf
b}^2}{{\sf c}^4\overline{\sf c}^2}+8\,\frac{P{\sf b}^2\overline{\sf
b}}{{\sf c}^4\overline{\sf c}^3}+ 9i\,\frac{{{\sf b}^2\overline{\sf
b}^2}}{{\sf c}^4\overline{\sf c}^4}+
\\
&+\bigg(-\frac{\mathcal L(\overline P)}{{\sf c}\overline{\sf
c}}+2i\,\frac{P\sf b}{{\sf c}^2\overline{\sf c}}-2i\,\frac{\overline
P\,\overline{\sf b}}{{\sf c}\overline{\sf c}^2}-6\,\frac{{{\sf
b}\overline{\sf b}}}{{\sf c}^2\overline{\sf c}^2}\bigg)\,{\sf s}.
\endaligned
\end{equation}

After determining so the group parameter $\sf r$,
we have to re-compute $\gamma$ which can now be expressed as a combination
of the lifted coframe
$\rho,\zeta,\overline\zeta,\alpha,\beta,\widetilde{\alpha},\widetilde{\beta}$
independently of $d\sf r$, {\it cf.} \thetag{\ref{gamma0}}.
For this, first we need the expression of $d\sf r$, not
only of $\sf r$.

Differentiating $\sf r$ in \thetag{\ref{r}} gives:
\[
\footnotesize\aligned d{\sf r}&=-\frac{1}{3\,{\sf c}\overline{\sf
c}^2}\,d\mathcal L(\overline{\mathcal L}(\overline P))+\frac{\mathcal
L(\overline{\mathcal L}(\overline P))}{3\,{\sf c}^2\overline{\sf
c}^2}\,d{\sf c}+\frac{2\,\mathcal L(\overline{\mathcal L}(\overline
P))}{3\,{\sf c}\overline{\sf c}^3}d\overline{\sf c}+\frac{1}{2\,{\sf
c}\overline{\sf c}^2}d\overline{\mathcal L}(\overline{\mathcal
L}(P))-\frac{\overline{\mathcal L}(\overline{\mathcal L}(P))}{2\,{\sf
c}^2\overline{\sf c}^2}d{\sf c}-\frac{\overline{\mathcal
L}(\overline{\mathcal L}(P))}{{\sf c}\overline{\sf
c}^3}d\overline{\sf c}-
\\
&-\frac{i}{2}\frac{\overline{\mathcal L}(P)}{{\sf c}^2\overline{\sf
c}^2}d{\sf b}-\frac{i}{2}\frac{\sf b}{{\sf c}^2\overline{\sf
c}^2}d\overline{\mathcal L}(P)+i\,\frac{\overline{\mathcal L}(P)\sf
b}{{\sf c}^3\overline{\sf c}^2}d{\sf c}+i\,\frac{\overline{\mathcal
L}(P)\sf b}{{\sf c}^2\overline{\sf c}^3}d\overline{\sf
c}-\frac{\overline P}{6\,{\sf c}\overline{\sf
c}^2}d\overline{\mathcal L}(P)-\frac{\overline{\mathcal
L}(P)}{6\,{\sf c}\overline{\sf c}^2}d\overline P+\frac{\overline
P\overline{\mathcal L}(P)}{6\,{\sf c}^2\overline{\sf c}^2}d{\sf
c}+\frac{\overline P\overline{\mathcal L}(P)}{3\,{\sf c}\overline
{\sf c}^3}d\overline{\sf c}+
\\
&+\frac{{\sf b}\overline{\sf b}}{{\sf c}^2\overline{\sf
c}^3}d\overline P+\frac{\overline P\sf b}{{\sf c}^2\overline{\sf
c}^3}d\overline{\sf b}+\frac{\overline P\,\overline{\sf b}}{{\sf
c}^2\overline{\sf c}^3}d{\sf b}-2\,\frac{\overline P{\sf
b}\overline{\sf b}}{{\sf c}^3\overline{\sf c}^3}d{\sf
c}-3\,\frac{\overline P{\sf b}\overline{\sf b}}{{\sf
c}^2\overline{\sf c}^4}d\overline{\sf c}-2i\,\frac{{\sf
b}\overline{\sf b}}{{\sf c}^3\overline{\sf c}^3}d{\sf
b}-i\,\frac{{\sf b}^2}{{\sf c}^3\overline{\sf c}^3}d\overline{\sf
b}+3i\,\frac{{\sf b}^2\overline{\sf b}}{{\sf c}^4\overline{\sf
c}^3}d{\sf c}+3i\,\frac{{\sf b}^2\overline{\sf b}}{{\sf
c}^3\overline{\sf c}^4},
\endaligned
\]
in which, similarly
to the expressions \thetag{\ref{dP}} and \thetag{\ref{dLP}}, one has to
replace the differentials:
\[
\footnotesize\aligned
 d\mathcal L(\overline{\mathcal L}(\overline P))&=\bigg(\frac{1}{\sf c}\mathcal
L(\mathcal L(\overline{\mathcal L}(\overline
P)))\bigg)\cdot\zeta+\bigg(\frac{1}{\overline{\sf
c}}\overline{\mathcal L}(\mathcal L(\overline{\mathcal L}(\overline
P)))\bigg)\cdot\overline\zeta+
\\
&+\bigg(-\frac{\sf b}{{\sf c}^2\overline{\sf c}}\mathcal L(\mathcal
L(\overline{\mathcal L}(\overline P)))-\frac{\overline{\sf b}}{{\sf
c}\overline{\sf c}^2}\overline{\mathcal L}(\mathcal
L(\overline{\mathcal L}(\overline P)))+\frac{1}{{\sf c}\overline{\sf
c}}\mathcal T(\mathcal L(\overline{\mathcal L}(\overline
P)))\bigg)\cdot\rho,
\\
d\overline{\mathcal L}(\overline{\mathcal L}(P))&=\bigg(\frac{1}{\sf
c}\mathcal L(\overline{\mathcal L}(\overline{\mathcal
L}(P)))\bigg)\cdot\zeta+\bigg(\frac{1}{\overline{\sf
c}}\overline{\mathcal L}(\overline{\mathcal L}(\overline{\mathcal
L}(P)))\bigg)\cdot\overline\zeta+
\\
&+\bigg(-\frac{\sf b}{{\sf c}^2\overline{\sf c}}\mathcal
L(\overline{\mathcal L}(\overline{\mathcal
L}(P)))-\frac{\overline{\sf b}}{{\sf c}\overline{\sf
c}^2}\overline{\mathcal L}(\overline{\mathcal L}(\overline{\mathcal
L}(P)))+\frac{1}{{\sf c}\overline{\sf c}}\mathcal
T(\overline{\mathcal L}(\overline{\mathcal L}(P)))\bigg)\cdot\rho.
\endaligned
\]
Then thanks to the expressions \thetag{\ref{db-dc}}, one can
re-express $d\sf r$ in terms of the lifted coframe
$\rho,\zeta,\overline\zeta,\alpha,\beta,\widetilde{\alpha},
\widetilde{\beta}$. Because of the length of the result, we do not
present this intermediate computation here. After all, replacing
${\sf r}$ and $d{\sf r}$ in the Maurer-Cartan form $\gamma$ in
\thetag{\ref{gamma0}} re-shapes its expression under the form:
\begin{equation}
\label{gamma-2} \aligned
\gamma:=V_1\,\rho+V_2\,\zeta+V_3\,\overline\zeta,
\endaligned
\end{equation}
with three certain functions given by:
\[
\footnotesize\aligned V_1&:=-\frac{1}{3}\frac{{\mathcal T}({\mathcal
L}(\overline{\mathcal L}(\overline P)))}{{\sf c}^2\overline{\sf
c}^3}+\frac{\mathcal T(\overline{\mathcal L}(\overline {\mathcal
L}(P)))}{{\sf c}^2\overline{\sf c}^3}+\frac{1}{3}\frac{{\mathcal
L}({\mathcal L}(\overline{\mathcal L}(\overline P))){\sf b}}{{\sf
c}^3\overline{\sf c}^3}+\frac{1}{3}\frac{\overline{\mathcal
L}({\mathcal L}(\overline{\mathcal L}(\overline P)))\overline{\sf
b}}{{\sf c}^2\overline{\sf c}^4}-\frac{1}{2}\frac{\overline{\sf
b}\overline{\mathcal L}(\overline{\mathcal L}({\mathcal L}(\overline
P)))}{{\sf c}^2\overline{\sf
c}^4}-\frac{1}{2}\frac{\overline{\mathcal L}({\mathcal L}({\mathcal
L}(\overline P)))\sf b}{{\sf c}^3\overline{\sf c}^3}+
\\
&+\frac{i}{6}\frac{\overline P\overline{\mathcal L}({\mathcal
L}({\mathcal L}(\overline P)))}{{\sf c}^2\overline{\sf
c}^3}-\frac{i}{6}\frac{\overline P{\mathcal L}(\overline{\mathcal
L}({\mathcal L}(\overline P)))}{{\sf c}^2\overline{\sf
c}^3}-3i\,\frac{{\sf b}^2\overline{\sf b}\sf s}{{\sf
c}^3\overline{\sf c}^3}-\frac{i}{3}\frac{\overline{\mathcal
L}({\mathcal L}(P)){\sf b}^2}{{\sf c}^4\overline{\sf c}^3}-
\frac{5i}{2}\frac{\overline{\mathcal L}({\mathcal L}(\overline
P)){\sf b}\overline{\sf b}}{{\sf c}^3\overline{\sf c}^4}-
\frac{{\mathcal L}(\overline{\mathcal L}(\overline P))\sf s}{{\sf
c}\overline{\sf c}^2}+\frac{2}{3}\frac{{\mathcal
L}(\overline{\mathcal L}(\overline P))P\sf b}{{\sf c}^3\overline{\sf
c}^3}-
\\
&-\frac{1}{3}\frac{{\mathcal L}(\overline{\mathcal L}(\overline
P))\overline P\overline{\sf b}}{{\sf c}^2\overline{\sf
c}^4}+\frac{3}{2}\frac{\overline{\mathcal L}({\mathcal L}(\overline
P))\sf s}{{\sf c}\overline{\sf c}^2}- \frac{\overline{\mathcal
L}({\mathcal L}(\overline P))P\sf b}{{\sf c}^3\overline{\sf
c}^3}+\frac{2}{3}\frac{\overline{\mathcal L}({\mathcal L}(\overline
P))\overline P\overline{\sf b}}{{\sf c}^2\overline{\sf
c}^4}-\frac{1}{3}\frac{{\mathcal L}({\mathcal L}(\overline
P))\overline P\sf b}{{\sf c}^3\overline{\sf
c}^3}+\frac{1}{3}\frac{\overline{\mathcal L}({\mathcal
L}(P))\overline P\sf b}{{\sf c}^3\overline{\sf
c}^3}+i\,\frac{{\mathcal L}({\mathcal L}(\overline P)){\sf b}^2}{{\sf
c}^4\overline{\sf c}^3}+
\\
&+\frac{7i}{3}\frac{{\mathcal L}(\overline{\mathcal L}(\overline
P)){\sf b}\overline{\sf b}}{{\sf c}^3\overline{\sf
c}^4}-\frac{i}{12}\frac{\overline{\mathcal L}({\mathcal L}(\overline
P)){\mathcal L}(\overline P)}{{\sf c}^2\overline{\sf c}^3}
-\frac{3i}{2}\frac{\mathcal L(\overline P)\sf bs}{{\sf
c}^2\overline{\sf c}^2}-5\,\frac{{\sf b}^3\overline{\sf b}^2}{{\sf
c}^5\overline{\sf c}^5}-\frac{\mathcal L(\overline P)\overline P\sf
s}{{\sf c}\overline{\sf c}^2}-\frac{1}{6}\frac{{\mathcal L}(\overline
P)\overline{\mathcal L}(\overline P)\overline{\sf b}}{{\sf
c}^2\overline{\sf c}^4}+\frac{1}{2}\frac{{\mathcal L}(\overline
P)P\overline P\sf b}{{\sf c}^3\overline{\sf c}^3}-
\\
&-\frac{1}{6}\frac{{\mathcal L}(\overline P)\overline
P^2\overline{\sf b}}{{\sf c}^2\overline{\sf c}^4}-
\frac{\overline{\mathcal L}(\overline P){\sf b}{\sf b}^2}{{\sf
c}^3\overline{\sf c}^5}-4\,\frac{{\mathcal L}(\overline P){\sf
b}^2\overline{\sf b}}{{\sf c}^4\overline{\sf c}^4}+3\,\frac{\overline
P{\sf b}\overline{\sf b}\sf s}{{\sf c}^2\overline{\sf
c}^3}+\frac{\overline P^2{\sf b}\overline{\sf b}^2}{{\sf
c}^3\overline{\sf c}^5}-\frac{1}{12}\frac{\mathcal L(\overline
P)\mathcal L(\overline P)\sf b}{{\sf c}^3\overline{\sf
c}^3}-3\,\frac{P\overline P{\sf b}^2\overline{\sf b}}{{\sf
c}^4\overline{\sf c}^4}+3i\,\frac{P{\sf b}^3\overline{\sf b}}{{\sf
c}^5\overline{\sf c}^4}+ \frac{5i}{6}\frac{P{\mathcal L}(\overline
P){\sf b}^2}{{\sf c}^4\overline{\sf c}^3}-
\\
&-6i\,\frac{\overline P{\sf b}^2\overline{\sf b}^2}{{\sf
c}^4\overline{\sf c}^5}+\frac{i}{12}\frac{{\mathcal L}(\overline
P){\mathcal L}(\overline P)\overline P}{{\sf c}^2\overline{\sf
c}^3}-\frac{5i}{6}\frac{{\mathcal L}(\overline P)\overline P{\sf
b}\overline{\sf b}}{{\sf c}^3\overline{\sf c}^4},
\\
V_2&:=\frac{1}{2}\frac{{\mathcal L}(\overline{\mathcal L}({\mathcal
L}(\overline P)))}{{\sf c}^2\overline{\sf c}^2}-
\frac{1}{3}\frac{{\mathcal L}({\mathcal L}(\overline{\mathcal
L}(\overline P)))}{{\sf c}^2\overline{\sf
c}^2}-\frac{P\overline{\mathcal L}(\mathcal L(\overline P))}{{\sf
c}^2\overline{\sf c}^2}+\frac{2}{3}\frac{P\mathcal
L(\overline{\mathcal L}(\overline P))}{{\sf c}^2\overline{\sf c}^2}-
\frac{2i}{3}\frac{\overline{\mathcal L}({\mathcal L}(P))\sf b}{{\sf
c}^3\overline{\sf c}^2}-\frac{1}{6}\frac{{\mathcal L}({\mathcal
L}(\overline P))\overline P}{{\sf c}^2\overline{\sf
c}^2}-i\,\frac{\overline{\mathcal L}({\mathcal L}(\overline
P))\overline{\sf b}}{{\sf c}^2\overline{\sf c}^3}+
\\
&+\frac{2i}{3}\frac{{\mathcal L}(\overline{\mathcal L}(\overline
P))\overline{\sf b}}{{\sf c}^2\overline{\sf c}^3} {\sf
s}^2-\frac{\mathcal L(P){\sf b}^2}{{\sf c}^4\overline{\sf
c}^2}+\frac{1}{3}\frac{P\mathcal L(\overline P)\overline P}{{\sf
c}^2\overline{\sf c}^2}+\frac{i}{3}\frac{\mathcal L(\overline
P)\overline P\overline{\sf b}}{{\sf c}^2\overline{\sf
c}^3}+\frac{2i}{3}\frac{{\mathcal L}(\overline P)P\sf b}{{\sf
c}^3\overline{\sf c}^2}-\frac{1}{6}\frac{{\mathcal L}(\overline
P){\mathcal L}(\overline P)}{{\sf c}^2\overline{\sf
c}^2}-2i\,\frac{P{\sf b}^2\overline{\sf b}}{{\sf c}^4\overline{\sf
c}^3}+2\,\frac{{\sf b}^2\overline{\sf b}^2}{{\sf c}^4\overline{\sf
c}^4},
\\
  V_3&:=\frac{1}{2}\frac{\overline{\mathcal L}(\overline{\mathcal
L}({\mathcal L}(\overline P)))}{{\sf c}\overline{\sf
c}^3}-\frac{1}{3}\frac{\overline{\mathcal L}({\mathcal
L}(\overline{\mathcal L}(\overline P)))}{{\sf c}\overline{\sf
c}^3}+\frac{2}{3}\frac{{\mathcal L}(\overline{\mathcal L}(\overline
P))\overline P}{{\sf c}\overline{\sf
c}^3}-\frac{7}{6}\frac{\overline{\mathcal L}({\mathcal L}(\overline
P))\overline P}{{\sf c}\overline{\sf c}^3}
-\frac{1}{6}\frac{{\mathcal L}(\overline P)\overline{\mathcal
L}(\overline P)}{{\sf c}\overline{\sf
c}^3}+\frac{1}{3}\frac{{\mathcal L}(\overline P)\overline P^2}{{\sf
c}\overline{\sf c}^3}.
\endaligned
\]
One should notice that $V_2$ depends on the group
parameter $\sf s$, while $V_1$ and $V_3$ do not.

Now, substituting this new expression of $\gamma$ into the lastly achieved
structure equation \thetag{\ref{modified-structure-equation-2}},
changes it into the form (remind that $W_2$ vanishes after
determining $\sf r$):
\begin{eqnarray}
\label{modified-structure-equation-3}
\begin{array}{ll}
d\rho=\alpha\wedge\rho+\widetilde{\alpha}\wedge\rho+i\,\zeta\wedge\overline\zeta,
&
\\
d\zeta=\beta\wedge\rho+\alpha\wedge\zeta, &
\\
d\overline\zeta=\widetilde{\beta}\wedge\rho+\widetilde{\alpha}\wedge\overline\zeta,
&
\\
d\alpha=\delta\wedge\rho+2\,i\,\zeta\wedge\overline{\beta}+i\,\overline{\zeta}\wedge\beta,
&
\\
d\beta=\delta\wedge\zeta+\beta\wedge\overline{\alpha}+V_2\,\zeta\wedge\rho+V_3\,\overline\zeta\wedge\rho
 &
 \\
\ \ \ \ \ \
=\big(\delta-V_2\,\rho\big)\wedge\zeta+\beta\wedge\overline{\alpha}+V_3\,\overline\zeta\wedge\rho,
 &
\\
d\widetilde{\alpha}=\delta\wedge\rho-2\,i\,\overline\zeta\wedge{\beta}-i\,{\zeta}\wedge\widetilde{\beta},
&
  \\
d\widetilde{\beta}=\delta\wedge\overline\zeta+\widetilde{\beta}\wedge{\alpha}+i\,W_1\,\rho\wedge\overline\zeta+\overline
V_3\,\zeta\wedge\rho+\overline V_2\,\overline\zeta\wedge\rho &
\\
\ \ \ \ \ \ =\big(\delta+i\,W_1\,\rho-\overline
V_2\,\rho\big)\wedge\overline\zeta+\widetilde{\beta}\wedge{\alpha}+\overline
V_3\,\zeta\wedge\rho.
 &
\end{array}
\end{eqnarray}

At present, we have just one group parameter $\sf s$.
The
complete absorption will be rigorously possible only if
the seemingly {\it implausible} identity:
\[
V_2=-iW_1+\overline V_2,
\]
would be satisfied, because it would enable us to modify-rename:
\[
\aligned
\delta
:=
&\,\delta-V_2\rho
\\
=
&\,
\delta+
\big(i\,W_1-\overline V_2\big)\,\rho
\endaligned
\]
such a substitution for $\delta$ having no effect on the preceding
wedge product $\delta\wedge\rho$ in $d\alpha$ and
$d\widetilde{\alpha}$.

\noindent
 We claim that the desired identity
holds. In fact after simplification, we obtain:
\begin{equation}
\label{W1-V2} \footnotesize\aligned \overline
V_2-iW_1-V_2&=\frac{1}{3\,{\sf c}^2\overline{\sf
c}^2}\bigg(-3\,{\mathcal L}(\overline{\mathcal L}({\mathcal
L}(\overline P)))+3\,\overline{\mathcal L}({\mathcal L}({\mathcal
L}(\overline P)))+{\mathcal L}({\mathcal L}(\overline{\mathcal
L}(\overline P)))-\overline{\mathcal L}(\overline{\mathcal
L}({\mathcal L}(P)))+
\\
& \ \ \ \ \ \ \ \ \ \ \ \ \ \ \ \ \ \ \ \ +P\overline{\mathcal
L}({\mathcal L}(\overline P))-P{\mathcal L}(\overline{\mathcal
L}(\overline P))-\overline P{\mathcal L}({\mathcal L}(\overline
P))+\overline P\overline{\mathcal L}({\mathcal L}(P))\bigg).
\endaligned
\end{equation}
Serendipitously, this imaginary expression is much simplified and it
does not include the group parameter $\sf s$. To show that it
vanishes identically, we need the following result:

\begin{Lemma}
\label{iterated}
(\cite{Merker-Sabzevari-CEJM}, Proposition 6.1)
\label{hhh} Let $H_1$ and $H_2$ be two vector fields on a manifolds
$M$ satisfying:
\[
\aligned &[H_1,[H_1,H_2]]=\Phi_1[H_1,H_2], \ \ \ \
[H_2,[H_1,H_2]]=\Phi_2[H_1,H_2],
\endaligned
\]
for some two certain functions $\Phi_1$ and $\Phi_2$. Then the
following four identities involving
third-order derivatives are satisfied:
\[
\small \aligned 0 & \overset{\rm I}{\equiv} -\,H_1(H_2(H_1(\Phi_2)))
+ 2\,H_2(H_1(H_1(\Phi_2))) - H_2(H_2(H_1(\Phi_1)))
-\,\Phi_2\,H_1(H_2(\Phi_1)) + \Phi_2\,H_2(H_1(\Phi_1)),
\endaligned
\]
\[
\small \aligned 0 & \overset{\rm II}{\equiv} -\,H_2(H_1(H_1(\Phi_2)))
+ 2\,H_1(H_2(H_1(\Phi_2))) - H_1(H_1(H_2(\Phi_2)))
-\,\Phi_1\,H_2(H_1(\Phi_2)) + \Phi_1\,H_1(H_2(\Phi_2)),
\endaligned
\]
\[
\small \aligned 0 & \overset{\rm III}{\equiv}
-\,H_1(H_1(H_1(\Phi_2))) + 2\,H_1(H_2(H_1(\Phi_1))) -
H_2(H_1(H_1(\Phi_1))) + \Phi_1\,H_1(H_1(\Phi_2)) -
\Phi_1\,H_2(H_1(\Phi_1)),
\endaligned
\]
\[
\small \aligned 0 & \overset{\rm IV}{\equiv} H_2(H_2(H_1(\Phi_2))) -
2\,H_2(H_1(H_2(\Phi_2))) + H_1(H_2(H_2(\Phi_2))) -
\,\Phi_2\,H_2(H_1(\Phi_2)) + \Phi_2\,H_1(H_2(\Phi_2)). \qed
\endaligned
\]
\end{Lemma}

\begin{Corollary}
The above expression \thetag{\ref{W1-V2}} of $\overline V_2-iW_1-V_2$
in fact vanishes identically.
\end{Corollary}

\proof Subtracting the equation II from I gives:
\[
\small\aligned 0&\equiv
3\,H_2(H_1(H_1(\Phi_2)))-3\,H_1(H_2(H_1(\Phi_2)))-H_2(H_2(H_1(\Phi_1)))+H_1(H_1(H_2(\Phi_2)))-
\\
&\ \ \
-\Phi_2\,H_1(H_2(\Phi_1))+\Phi_2\,H_2(H_1(\Phi_1))+\Phi_1\,H_2(H_1(\Phi_2))-\Phi_1\,H_1(H_2(\Phi_2)).
\endaligned
\]
Now, it suffices to put $\Phi_1:=P, \Phi_2:=\overline P$ and
$H_1:=\mathcal L, H_2:=\overline{\mathcal L}$ into the above
equation, taking account of the reality condition $\mathcal
L(\overline P)=\overline{\mathcal L}(P)$.
\endproof

Consequently, the equality
$\delta-V_2\rho=\delta+iW_1\rho-\overline V_2\rho$ permits us
to apply the substitution $\delta\mapsto\delta-V_2\rho$.
After renaming the single torsion coefficient $V_3$ as
$\frak I$, the structure equations
\thetag{\ref{modified-structure-equation-3}} received
the much simplified form:
\begin{eqnarray}
\label{modified-structure-equation-4}
\begin{array}{ll}
d\rho=\alpha\wedge\rho+\widetilde{\alpha}\wedge\rho+i\,\zeta\wedge\overline\zeta,
&
\\
d\zeta=\beta\wedge\rho+\alpha\wedge\zeta, &
\\
d\overline\zeta=\widetilde{\beta}\wedge\rho+\widetilde{\alpha}\wedge\overline\zeta,
&
\\
d\alpha=\delta\wedge\rho+2\,i\,\zeta\wedge\overline{\beta}+i\,\overline{\zeta}\wedge\beta,
&
\\
d\beta=\delta\wedge\zeta+\beta\wedge\overline{\alpha}+{\frak
I}\,\overline\zeta\wedge\rho, &
\\
d\widetilde{\alpha}=\delta\wedge\rho-2\,i\,\overline\zeta\wedge{\beta}-i\,{\zeta}\wedge\widetilde{\beta},
 &
 \\
d\widetilde{\beta}=\delta\wedge\overline\zeta+\widetilde{\beta}\wedge{\alpha}+\overline{\frak
I}\,\zeta\wedge\rho,
 &
\end{array}
\end{eqnarray}
with the single (modified) Maurer-Cartan form $\delta$
(after simplification):
\begin{equation}
\label{delta-2}\footnotesize\aligned \delta&=d{\sf s}+
\\
&+\bigg(-{\sf s}^2+\frac{1}{3}\frac{{\mathcal L}({\mathcal
L}(\overline{\mathcal L}(\overline P)))}{{\sf c}^2\overline{\sf
c}^2}-\frac{1}{2}\frac{{\mathcal L}(\overline{\mathcal L}({\mathcal
L}(\overline P)))}{{\sf c}^2\overline{\sf
c}^2}-\frac{2}{3}\frac{P{\mathcal L}(\overline{\mathcal L}(\overline
P))}{{\sf c}^2\overline{\sf
c}^2}+\frac{2i}{3}\frac{\overline{\mathcal L}({\mathcal L}(P))\sf
b}{{\sf c}^3\overline{\sf c}^2}+\frac{\overline{\mathcal L}({\mathcal
L}(\overline P))P}{{\sf c}^2\overline{\sf c}^2}+
\\
&+i\,\frac{\overline{\mathcal L}({\mathcal L}(\overline
P))\overline{\sf b}}{{\sf c}^2\overline{\sf
c}^3}+\frac{1}{6}\frac{{\mathcal L}({\mathcal L}(\overline
P))\overline P}{{\sf c}^2\overline{\sf
c}^2}-\frac{2i}{3}\frac{{\mathcal L}(\overline{\mathcal L}(\overline
P))\overline{\sf b}}{{\sf c}^2\overline{\sf
c}^3}-\frac{2i}{3}\frac{{\mathcal L}(\overline P)P\sf b}{{\sf
c}^3\overline{\sf c}^2}+\frac{{\mathcal L}(P){\sf b}^2}{{\sf
c}^4\overline{\sf c}^2}-
\\
&-\frac{1}{3}\frac{{\mathcal L}(\overline P)P\overline P}{{\sf
c}^2\overline{\sf c}^2}-\frac{i}{3}\frac{{\mathcal L}(\overline
P)\overline P\,\overline{\sf b}}{{\sf c}^2\overline{\sf
c}^3}+\frac{1}{6}\frac{{\mathcal L}(\overline P){\mathcal
L}(\overline P)}{{\sf c}^2\overline{\sf c}^2}-2\frac{{\sf
b}^2\overline{\sf b}^2}{{\sf c}^4\overline{\sf c}^4}+2i\,\frac{P{\sf
b}^2\overline{\sf b}}{{\sf c}^4\overline{\sf c}^3}\bigg)\cdot\rho+
\\
&+\bigg(\frac{P\sf s}{\sf c}+2i\,\frac{{\sf s}\overline{\sf b}}{{\sf
c}\overline{\sf c}}-\frac{i}{3}\frac{\overline{\mathcal L}({\mathcal
L}(P))}{{\sf c}^2\overline{\sf c}}+\frac{i}{3}\frac{{\mathcal
L}(\overline P)P}{{\sf c}^2\overline{\sf c}}-\frac{{\mathcal L}(P)\sf
b}{{\sf c}^3\overline{\sf c}}+2\,\frac{{\sf b}\overline{\sf
b}^2}{{\sf c}^3\overline{\sf c}^3}-2i\,\frac{P{\sf b}\overline{\sf
b}}{{\sf c}^3\overline{\sf c}^2}\bigg)\cdot\zeta+
\\
 &+\bigg(i\,\frac{\sf bs}{{\sf c}\overline{\sf
c}}-\frac{i}{2}\frac{\overline{\mathcal L}({\mathcal L}(\overline
P))}{{\sf c}\overline{\sf c}^2}+\frac{i}{3}\frac{{\mathcal
L}(\overline{\mathcal L}(\overline P))}{{\sf c}\overline{\sf c}^2}+
\frac{i}{6}\frac{{\mathcal L}(\overline P)\overline P}{{\sf
c}\overline{\sf c}^2}+2\,\frac{{\sf b}^2\overline{\sf b}}{{\sf
c}^3\overline{\sf c}^3}-i\,\frac{P{\sf b}^2}{{\sf c}^3\overline{\sf
c}^2}\bigg)\cdot\overline\zeta+{\sf s}\,\alpha-\bigg(\frac{P}{\sf
c}+2i\,\frac{\overline{\sf b}}{{\sf c}\overline{\sf
c}}\bigg)\cdot\beta+{\sf s}\,\widetilde{\alpha}-i\,\frac{\sf b}{{\sf
c}\overline{\sf c}}\,\widetilde{\beta}.
\endaligned
\end{equation}
As mentioned before, $\frak I$ is independent
of the only remaining group parameter $\sf s$, hence it is
impossible to normalize it. Consequently, this torsion coefficient is
actually an {\it essential invariant} of the problem.

\subsection{Second prolongation}

In the situation that we have still one undetermined group parameter
$\sf s$ without the possibility of normalizing the single essential torsion
coefficient $\frak I$, we have to prolong the latest structure
equations \thetag{\ref{modified-structure-equation-4}} by adding the
group parameter $\sf s$ to the set of base variables $z,\overline
z,u,{\sf b},\overline{\sf b},{\sf c},\overline{\sf c}$ and adding the
1-form $\delta$ to the coframe
$\{\rho,\zeta,\overline\zeta,\alpha,\widetilde{\alpha},\beta,\widetilde{\beta}\}$.
Before starting this step, let us present the following result:

\begin{Lemma}
The above modified 1-form $\delta$ is the unique one which enjoys the
structure equations \thetag{\ref{modified-structure-equation-4}}.
\end{Lemma}

\proof Assume that $\delta$ and $\delta'$ are two forms satisfying the
structure equations, simultaneously. A subtraction
immediately gives:
\[
\aligned 0\equiv(\delta-\delta')\wedge\rho,
\ \ \ \ \ \ \ \ \ \ \ \ \ \ \ \
0\equiv(\delta-\delta')\wedge\zeta,
\endaligned
\]
which according to Cartan's lemma implies that $\delta-\delta'$ must
be a combination of only $\rho$ and of only $\zeta$, which
clearly implies $\delta - \delta' = 0$.
\endproof

This shows that we do not encounter any new
(prolonged) group parameter while starting the next prolongation.
In other words, the prolonged structure group will be
automatically reduced
to an $e$-{\it structure}. Hence
it only remains to compute $d \delta$.

\begin{Proposition}
The exterior differentiation $d\delta$  has the form:
\begin{equation}
\label{ddelta-lem}
d\delta=\delta\wedge\alpha+\delta\wedge\widetilde{\alpha}+i\,\beta\wedge\widetilde{\beta}+{\frak
T}\,\rho\wedge\zeta+\overline{\frak T}\,\rho\wedge\overline\zeta,
\end{equation}
for a certain complex function $\frak T$.
\end{Proposition}

\proof Differentiating $d\alpha$ in the last structure equation
\thetag{\ref{modified-structure-equation-4}} gives:
\begin{equation*}\footnotesize
\aligned\footnotesize 0&\equiv
d\delta\wedge\rho-\delta\wedge\alpha\wedge\rho-\delta\wedge\overline{\alpha}\wedge\rho
\underline{-i\,\delta\wedge\zeta\wedge\overline{\zeta}}_{\scriptsize{\sf
a}}\underline{-2i\,\delta\wedge\overline{\zeta}\wedge\zeta}_{\scriptsize{\sf
a}}\underline{-2i\,\overline{\beta}\wedge\alpha\wedge\zeta}_{\scriptsize{\sf
b}}+
\\
&\underline{+2\,i\,\overline{\beta}\wedge\alpha\wedge\zeta}_{\scriptsize{\sf
b}}+\underline{2i\,\widetilde{\beta}\wedge\beta\wedge\rho}_{\scriptsize{\sf
c}}\underline{-i\,\delta\wedge\zeta\wedge\overline\zeta}_{\scriptsize{\sf
a}}\underline{-i\,\beta\wedge\overline{\alpha}\wedge\overline\zeta}_{\scriptsize{\sf
d}}\underline{+i\,\beta\wedge\overline{\alpha}\wedge\overline{\zeta}}_{\scriptsize{\sf
d}}+\underline{i\,\beta\wedge\overline{\beta}\wedge\rho}_{\scriptsize{\sf
c}},
\endaligned
\end{equation*}
in which the underlined terms can be simplified together and bring
the following simple equality:
\begin{eqnarray}
\label{dphi}
(d\delta-\delta\wedge\alpha-\delta\wedge\overline{\alpha}-i\,\beta\wedge\overline{\beta})\wedge\rho\equiv
0.
\end{eqnarray}
On the other hand, from differentiating $d\beta$ and
$d\overline{\beta}$ we also find:
\begin{equation}
\label{dpsi} \aligned
&(d\delta-\delta\wedge\alpha-\delta\wedge\overline{\alpha}-i\,\beta\wedge\overline{\beta})\wedge\zeta+
(\underbrace{d {\frak I}\wedge\overline{\zeta}-3\,{\frak
I}\,\overline{\alpha}\wedge\overline{\zeta}+{\frak I}\,\alpha\wedge
\overline{\zeta}}_{\Gamma})\wedge\rho\equiv 0,
\\
&(d\delta-\delta\wedge\alpha-\delta\wedge\overline{\alpha}-i\,\beta\wedge\overline{\beta})\wedge\overline{\zeta}+
(\underbrace{d\overline{{\frak I}}\wedge\zeta-3\,\overline{{\frak
I}}\,\alpha\wedge\zeta+\overline{{\frak I}}\,\overline{\alpha}\wedge
\zeta}_{\overline{\Gamma}})\wedge\rho\equiv 0,
\endaligned
\end{equation}
after a slight simplification.
Now, applying the Cartan's Lemma \ref{Cartan's Lemma} to the equality
\thetag{\ref{dphi}} gives:
\begin{eqnarray}
\label{eta}
d\delta=\delta\wedge\alpha+\delta\wedge\overline{\alpha}+i\,\beta\wedge\overline{\beta}+\xi\wedge\rho,
\end{eqnarray}
for some 1-form $\xi$. Putting then this
expression of $d\delta$ into \thetag{\ref{dpsi}}
brings:
\begin{eqnarray}\footnotesize
\label{Gamma-xi2}
(\xi\wedge\zeta-\Gamma)\wedge\rho&=&0,
\\
\nonumber
(\xi\wedge\overline{\zeta}-\overline{\Gamma})\wedge\rho&=&0.
\end{eqnarray}
Applying again the Cartan's Lemma to the first equation, we get:
\[
\xi\wedge\zeta-\Gamma=\mathcal A\wedge\rho,
\]
for some 1-form $\mathcal A$,
or equivalently:
\[
\xi\wedge\zeta-(d{\frak I}-3\,{\frak I}\overline{\alpha}+{\frak
I}\alpha)\wedge\overline{\zeta}-\mathcal A\wedge\rho=0.
\]
Applying the Cartan's Lemma, this time
to the last equality, we obtain:
\begin{eqnarray}
\label{xi} \xi=A_1\zeta+A_2\overline{\zeta}+A_3\rho,
\end{eqnarray}
for some certain functions $A_1,A_2$, $A_3$. Subtracting the
conjugation of the second equation in \thetag{\ref{Gamma-xi2}} from
the first one also gives:
\[
(\xi\wedge\zeta-\overline{\xi}\wedge\zeta)\wedge\rho\equiv 0,
\]
and hence there is a 1-form $\mathcal C$ with:
\[
(\xi-\overline{\xi})\wedge\zeta+{\mathcal C}\wedge\rho\equiv 0.
\]
We apply again the Cartan's lemma and this time we obtain the
following equation for two certain complex functions $B_1$ and $B_2$:
\begin{eqnarray}
\label{xi2}
 \xi-\overline{\xi}=B_1\zeta+B_2\rho.
\end{eqnarray}
The left-hand side of this equality is imaginary and hence the coefficient of
$\zeta$ must vanish:
$B_1=0$. On the other hand, according to \thetag{\ref{xi}}
we have:
\[
\xi-\overline{\xi}=(A_1-\overline{A}_2)\zeta+(A_2-\overline{A}_1)\overline{\zeta}+(A_3-\overline{A}_3)\rho.
\]
Comparing this equation with \thetag{\ref{xi2}}
then immediately implies that $A_2=\overline{A}_1$.
Hence, denoting $-A_1$ by ${\frak T}$ gives the following expression
for the 2-form $\xi\wedge\rho$ according to \thetag{\ref{xi}}:
\[
\xi\wedge\rho={\frak T}\rho\wedge\zeta+\overline{{\frak
T}}\,\rho\wedge\overline{\zeta}.
\]
To complete the proof, it is now enough to put the above expression
into \thetag{\ref{eta}}.
\endproof

Consequently we will have the following (prolonged) structure
equations after adding the differentiation of the new lifted 1-form
$\delta$ to the previous ones:
\begin{eqnarray}
\label{modified-structure-equation-final}
\begin{array}{ll}
d\rho=\alpha\wedge\rho+\widetilde{\alpha}\wedge\rho+i\,\zeta\wedge\overline\zeta,
&
\\
d\zeta=\beta\wedge\rho+\alpha\wedge\zeta, &
\\
d\overline\zeta=\widetilde{\beta}\wedge\rho+\widetilde{\alpha}\wedge\overline\zeta,
&
\\
d\alpha=\delta\wedge\rho+2\,i\,\zeta\wedge\overline{\beta}+i\,\overline{\zeta}\wedge\beta,
&
\\
d\beta=\delta\wedge\zeta+\beta\wedge\overline{\alpha}+{\frak
I}\,\overline\zeta\wedge\rho,
 &
\\
d\widetilde{\alpha}=\delta\wedge\rho-2\,i\,\overline\zeta\wedge{\beta}-i\,{\zeta}\wedge\widetilde{\beta},
&

 \\
d\widetilde{\beta}=\delta\wedge\overline\zeta+\widetilde{\beta}\wedge{\alpha}+\overline{\frak
I}\,\zeta\wedge\rho,
 &
 \\
 d\delta=\delta\wedge\alpha+\delta\wedge\widetilde{\alpha}+i\,\beta\wedge\widetilde{\beta}+{\frak
T}\,\rho\wedge\zeta+\overline{\frak T}\,\rho\wedge\overline\zeta.
&
\end{array}
\end{eqnarray}
These equations provide the final $e$-structure.

\smallskip

Our ultimate task is to find the expression of the new
coefficient $\frak T$. For this aim, we employ the same procedure as
that of finding the expression of $W$ in \thetag{\ref{W}}. At first,
we have to compute the exterior differentiatial of $\delta$ in
\thetag{\ref{delta-2}}. Unfortunately, this expression is extensive
(almost 2 pages long), hence we do not present it here.

Another much shorther path is to carefully compare this expression of
$d\delta$ to that
from~\thetag{\ref{modified-structure-equation-final}}.  Considering
the coefficient of $\rho\wedge\zeta$ reveals a compact expression for
$\frak T$, granted the four equations I--IV of Lemma \ref{iterated}
and their first order derivations with respect to the operators
$\mathcal L$ and $\overline{\mathcal L}$. Then one finds out that the
desired function $\frak T$ can be expressed in terms of the essential
invariant $\frak I$ as:
\begin{equation*}
\label{T} {\frak T}=\frac{1}{\overline{\sf
c}}\,\bigg(\overline{\mathcal L}(\overline{\frak I})-\overline
P\,\overline{\frak I}\bigg)-i\,\frac{\sf b}{{\sf c}\overline{\sf
c}}\overline{\frak I}.
\end{equation*}

Now, from standard features of the theory, we conclude:

\begin{Theorem}
The equivalence problem for strongly pseudoconvex
Levi-nondegenerate hypersurfaces $M^3\subset\mathbb C^2$ has a single
essential primary invariant:
\[\footnotesize\aligned
{\frak I}&=-\frac{1}{3}\frac{\overline{\mathcal L}({\mathcal
L}(\overline{\mathcal L}(\overline P)))}{{\sf c}\overline{\sf
c}^3}+\frac{2}{3}\frac{{\mathcal L}(\overline{\mathcal L}(\overline
P))\overline P}{{\sf c}\overline{\sf
c}^3}+\frac{1}{2}\frac{\overline{\mathcal L}(\overline{\mathcal
L}({\mathcal L}(\overline P)))}{{\sf c}\overline{\sf
c}^3}-\frac{7}{6}\frac{\overline{\mathcal L}({\mathcal L}(\overline
P))\overline P}{{\sf c}\overline{\sf c}^3}-
\\
&\ \ \ \ \ -\frac{1}{6}\frac{{\mathcal L}(\overline
P)\overline{\mathcal L}(\overline P)}{{\sf c}\overline{\sf
c}^3}+\frac{1}{3}\frac{{\mathcal L}(\overline P)\overline P^2}{{\sf
c}\overline{\sf c}^3},
\endaligned
\]
in which the fundamental function $P:=P(z,\overline z,u)$ expresses
explicitly in terms of the graphing function $\varphi$ as:
\[
P:=\frac{\ell_z-\ell A_u+A\ell_u}{\ell},
\]
where:
\[
A := \frac{i\,\varphi_z}{1-i\,\varphi_u} \ \ \ \ {\rm and} \ \ \ \
\ell := i\,\big( \overline{A}_z + A\,\overline{A}_u -
A_{\overline{z}} - \overline{A}\,A_u \big).
\]
In particular, this invariant vanishes when and
only when $M^3$ is biholomorphic to the model
Heisenberg sphere defined as the graph of the function:
\[
v=z\overline z.
\]

\proof
It is only necessary to observe that with the assumption
$\varphi(z,\overline z,u) := z\overline z$,
one immediately gets $P\equiv
0$, and hence $\frak I \equiv 0$.

Conversely, if $\mathfrak{I} = 0$, whence also $\mathfrak{ T} = 0$,
the constructed $e$-structure identifies with the Maurer-Cartan
equations of the real projective group, and one recovers the
Heisengerg sphere as the orbit of the origin under the action of this
group.
\endproof

\end{Theorem}

\section{A brief comparison to the
\\
Cartan-Tanaka geometry of real hypersurfaces $M^3\subset\mathbb C^2$}
\label{Cartan geometry}

We now turn to a brief discussion of Cartan geometry of the under
consideration real hypersurfaces $M^3\subset\mathbb C^2$  which is
much pertinent to their problem of equivalence. It helps us to
understand better the generally close relationship between the
equivalence problems and Cartan geometries. Here, we borrow the
results, notations and terminology from the recent paper
\cite{Merker-Sabzevari-CEJM} ({\it see} also \cite{Sharpe-1997}).

\begin{Definition}
\label{Cartanconnecdefini} Let $G$ be a Lie group with a closed
subgroup $H$, and let $\frak g$ and $\frak h$ be the corresponding
Lie algebras. A {\sl Cartan geometry of type $(G,H)$} on a manifold
$M$ is a principal $H$-bundle:
\[
\pi:\mathcal G\longrightarrow M
\]
together with a $\frak g$-valued $1$-form $\omega$, called the
corresponding {\sl Cartan connection}, on $\mathcal G$ subjected to
the following three conditions:

\smallskip\begin{itemize}

\item[\textbf{(i)}] $\omega_p:T_p\mathcal G\longrightarrow\frak g$ is
a linear isomorphism at every point $p\in\mathcal G$;

\smallskip \item[\textbf{(ii)}]
if $R_h(p):=ph$ is the right translation on $\mathcal G$ by any $h\in
H$, then:
\[
R^\ast_h\omega={\rm Ad}(h^{-1})\circ\omega;
\]

\smallskip\item[\textbf{(iii)}] $\omega(H^\dag)={\sf h}$ for every
${\sf h}\in\frak h$, where:
\[
H^\dag|_p := \textstyle{\frac{d}{dt}}\big|_0\big((R_{\exp(t \sf
h)}(p)\big)
\]
is the left-invariant vector field on $\mathcal G$ corresponding to
$\sf h$.
\end{itemize}

\noindent
\end{Definition}

Among Cartan geometries of type $(G,H)$, the most
symmetric one, called {\sl Klein geometry of type $(G,H)$}, arises
when $M=G/H$, when $\pi \colon G \rightarrow G/H$ is the projection
onto left-cosets, and when $\omega=\omega_{MC} \colon TG \to \frak g$
is the {\sl Maurer-Cartan form} on $G$.

In general, with a Cartan connection $\omega$ as above, if we
associate the vector field $\widehat{X}:=\omega^{-1}({\sf x})$ on
$\mathcal G$ to an arbitrary element ${\sf x}$ of $\frak g$, then the
infinitesimal version of condition \textbf{(ii)} reads as:
\[
[\widehat{X},\widehat{Y}] = \widehat{[{\sf x},{\sf y}]_{\frak g}},
\]
whenever {\sf y} belongs to $\frak h$. But in the special case of
Klein geometries, this equality holds moreover for any arbitrary
element $\sf y$ of $\frak g$. This difference motivates one to define
the {\sl curvature function}:
\[
\kappa \colon\ \ \ {\mathcal G} \longrightarrow {\rm Hom} \big(
\Lambda^2 (\frak g/\frak h),\frak g \big)
\]
associated to the Cartan connection $\omega$ by:
\[
\kappa_p({\sf x},{\sf y}) :=
 \omega_p\big([\widehat{X},
\widehat{Y}]\big)-[{\sf x},{\sf y}]_{\frak g} \ \ \ \ \ \ \ \ \ \ \ \
\ {\scriptstyle{(p\,\in\,\mathcal{G}, \ \ {\sf x},\,{\sf
y}\,\in\,\frak {g/h})}}.
\]
In a way, the curvature function measures how far a Cartan geometry
is from its corresponding Klein geometry. In particular, a Cartan
geometry is locally equivalent to its corresponding Klein geometry if
and only if its curvature function vanishes identically ({\em see}
\cite{Sharpe-1997}).

Now, let us return to the Levi-nondegenerate real hypersurfaces $M^3$
regarded as deformations of the Heisenberg sphere $\mathbb H^3$. In
\cite{Merker-Sabzevari-CEJM}, we built a regular
normal Cartan connection of type $(G,H)$ in which
$G$ is the projective group associated to the $8$-dimensional projective Lie
algebra:
\[
\frak g
:=
\frak{aut}(\mathbb H^3)={\sf Span}_{\mathbb R}({\sf
t,h}_1,{\sf h}_2,{\sf d,r},{\sf i}_1,{\sf i}_2,{\sf j})
\]
of
infinitesimal CR-automorphisms of $\mathbb H^3$ equipped with the
full commutator table:
\begin{center}
\footnotesize
\begin{tabular} [t] { l | l l l l l l l l }
\label{structure-g} &  ${\sf t}$ & ${\sf h}_1$ & ${\sf h}_2$ & ${\sf
d}$ & ${\sf r}$ & ${\sf i}_1$ & ${\sf i}_2$ & ${\sf j}$
\\
\hline ${\sf t}$ & $0$ &  $0$ & $0$ & $2\,{\sf t}$ & $0$ & ${\sf
h}_1$ & ${\sf h}_2$ & ${\sf d}$
\\
${\sf h}_1$ & $*$ & $0$ &  $4\,{\sf t}$ & ${\sf h}_1$ & ${\sf h}_2$ &
$6\,{\sf r}$ & $2\,{\sf d}$ & ${\sf i}_1$
\\
${\sf h}_2$ & $*$ & $*$ &  $0$ & ${\sf h}_2$ & $-{\sf h}_1$ &
$-2\,{\sf d}$ & $6\,{\sf r}$ & ${\sf i}_2$
\\
${\sf d}$ & $*$ & $*$ & $*$ &  $0$ & $0$ & ${\sf i}_1$ & ${\sf i}_2$
& $2\,{\sf j}$
\\
${\sf r}$ & $*$ & $*$ & $*$ &  $*$ & $0$ & $-{\sf i}_2$ & ${\sf i}_1$
& $0$
\\
${\sf i}_1$ & $*$ & $*$ & $*$ & $*$ & $*$ & $0$ & $4\,{\sf j}$ & $0$
\\
${\sf i}_2$ & $*$ & $*$ & $*$ & $*$ & $*$ & $*$ & $0$ & $0$
\\
${\sf j}$ & $*$ & $*$ & $*$ & $*$ & $*$ & $*$ & $*$ & $0$.
\end{tabular}
\end{center}
Moreover, $H$ is the
subgroup of $G$ associated to the 5-dimensional subalgebra:
\[
\frak
h:={\sf Span}_{\mathbb R}({\sf d,r},{\sf i}_1,{\sf i}_2,{\sf j}).
\]
The Lie algebra $\frak g$ is in fact graded, in the sense of Tanaka
\cite{Tanaka}:
\[
\frak g=\underbrace{\frak g_{-2}\oplus\frak g_{-1}}_{\frak
g_-}\oplus\underbrace{\frak g_{0}\oplus\frak g_1\oplus\frak
g_2}_{\frak h},
\]
with $\frak g_{-2}:={\sf Span}_{\mathbb R}({\sf t})$, with $\frak
g_{-1}:={\sf Span}_{\mathbb R}({\sf h}_1,{\sf h}_2)$, with $\frak
g_0:={\sf Span}_{\mathbb R}({\sf d,r})$, with $\frak g_1:={\sf
Span}_{\mathbb R}({\sf i}_1,{\sf i}_2)$ and with $\frak g_{2}:={\sf
Span}_{\mathbb R}({\sf j})$. Here $\frak g_-= \frak g \big/ \frak h$
is in fact the Levi-Tanaka symbol algebra of
any Levi nondegenerate $M^3 \subset \C^2$.

According to
this grading, the curvature function $\kappa$ decomposes into
homogeneous components:
\[
\kappa:=\kappa^{(0)}+\cdots+\kappa^{(5)}
\]
where $\kappa^{(s)}$
assigns to each pair $({\sf p}_{j_1},{\sf p}_{j_2})\in\Lambda^2\frak
g_-$, for ${\sf p}_{j_1}\in\frak g_{j_i}, j_i=-2,-1$, an element of
$\frak g_{j_1+j_2+s}$. It turns out that each curvature component
$\kappa^{(s)}$ can be formulated in the form:
\begin{equation}
\label{k-components} \aligned
\kappa^{(s)}=\sum_{s=j-(j_1+j_2)}\kappa^{p_{j_1}p_{j_2}}_{q_j} \,
{\sf p}_{j_1}^\ast\wedge {\sf p}^\ast_{j_2}\otimes {\sf q}_j,
\endaligned
\end{equation}
where $\kappa^{p_{j_1}p_{j_2}}_{q_j}(p)$ is the real-valued function
defined on an arbitrary point $p$ of $\mathcal G$ as the coefficient
of ${\sf q}_j$ in $\kappa(p)({\sf p}_{j_1},{\sf p}_{j_2})$, where
${\sf p}_{j_1}\in\frak g_{j_1},{\sf p}_{j_2}\in\frak g_{j_2}$, ${\sf
q}_j\in\frak g_j$ are some mentioned basis elements of $\frak g$, for
$j_1,j_2 = -2, -1$ and $j=-2, -1, 0, 1, 2$.

In fact, the process of construction the sought Cartan geometry in
\cite{Merker-Sabzevari-CEJM} has mainly consisted in annihilating as
many curvature components as possible, and finally we were able to
annihilate $\kappa^{ (0)}$ (easiest thing), $\kappa^{ (1)}$,
$\kappa^{ (2)}$ and $\kappa^{ (3)}$ by an appropriate progressive
building of $\omega$ which requires somewhat hard elimination
computations. Such computations have been done in the framework of
the powerful algorithm of Tanaka \cite{Tanaka} which involves some
modern concepts such as Lie algebras of infinitesimal
CR-automorphisms, Lie algebra cohomology, Tanaka prolongation and so
on. Finally we found out that (Proposition 7.3 and Theorem 7.4 of
\cite{Merker-Sabzevari-CEJM}):

\begin{Theorem}
The Cartan geometry associated to any
$\mathcal{ C}^6$-smooth Levi nondegenerate deformation
$M^3\subset\mathbb C^2$ of the Heisenberg sphere $\mathbb
H^3\subset\mathbb C^2$ has the curvature function:
\begin{equation}
\label{kappa}
 \aligned \kappa&=\kappa^{(4)}+\kappa^{(5)}=
\\
&= \kappa^{h_1t}_{i_1}{\sf h}_1^\ast\wedge{\sf t}^\ast\otimes {\sf
i}_1+\kappa^{h_1t}_{i_2}{\sf h}_1^\ast\wedge{\sf t}^\ast\otimes {\sf
i}_2+\kappa^{h_2t}_{i_1}{\sf h}_2^\ast\wedge{\sf t}^\ast\otimes {\sf
i}_1+
\\
&+ \kappa^{h_2t}_{i_2}{\sf h}_2^\ast\wedge{\sf t}^\ast\otimes {\sf
i}_2+\kappa^{h_1t}_j{\sf h}_1^\ast\wedge{\sf t}^\ast\otimes {\sf
j}+\kappa^{h_2t}_j{\sf h}_2^\ast\wedge{\sf t}^\ast\otimes {\sf j},
\endaligned
\end{equation}
with:
\[
\footnotesize \aligned \kappa_{i_1}^{h_1t} & =
-\,\mathbf{\Delta_1}\,c^4 - 2\,\mathbf{\Delta_4}\,c^3d -
2\,\mathbf{\Delta_4}\,cd^3 + \mathbf{\Delta_1}\,d^4,
\\
\kappa_{i_2}^{h_1t} & = -\,\mathbf{\Delta_4}\,c^4 +
2\,\mathbf{\Delta_1}\,c^3d + 2\,\mathbf{\Delta_1}\,cd^3 +
\mathbf{\Delta_4}\,d^4,
\\
\kappa_{i_1}^{h_2t} & = \kappa_{i_2}^{h_1t}, \ \ \ \ \ \
\kappa_{i_2}^{h_2t} = -\,\kappa_{i_1}^{h_1t},
\\
\kappa^{h_1t}_j & = \widehat{H}_1\big(\kappa^{h_2t}_{i_2}\big) -
\widehat{H}_2\big(\kappa^{h_1t}_{i_2}\big), \ \ \ \ \ \
\kappa^{h_2t}_j  = -\widehat{H}_1\big(\kappa^{h_2t}_{i_1}\big) +
\widehat{H}_2\big(\kappa^{h_1t}_{i_1}\big)
\endaligned
\]
and with the essential invariants, explicitly expressed in terms of
the defining function $\varphi$, as:
\begin{equation}
\label{Deltas} \footnotesize\aligned \mathbf{\Delta_1} & =
{\textstyle{\frac{1}{384}}} \Big[ H_1(H_1(H_1(\Phi_1))) -
H_2(H_2(H_2(\Phi_2))) + 11\,H_1(H_2(H_1(\Phi_2))) -
11\,H_2(H_1(H_2(\Phi_1)))+6\,\Phi_2\,H_2(H_1(\Phi_1))-
\\
& \ \ \ \ \ \ \ \ \ \ \ \ \  - 6\,\Phi_1\,H_1(H_2(\Phi_2)) -
3\,\Phi_2\,H_1(H_1(\Phi_2)) +
3\,\Phi_1\,H_2(H_2(\Phi_1))-\,3\,\Phi_1\,H_1(H_1(\Phi_1))+
3\,\Phi_2\,H_2(H_2(\Phi_2))-
\\
& \ \ \ \ \ \ \ \ \ \ \ \ \   - H_1(\Phi_1)\,H_1(\Phi_1) +
H_2(\Phi_2)\,H_2(\Phi_2) -\,2\,(\Phi_2)^2\,H_1(\Phi_1) +
2\,(\Phi_1)^2\,H_2(\Phi_2) - 2\,(\Phi_2)^2\,H_2(\Phi_2) +
2\,(\Phi_1)^2\,H_1(\Phi_1) \Big],
\endaligned
\end{equation}
\begin{equation*}
\footnotesize\aligned
 \mathbf{\Delta_4} & =
{\textstyle{\frac{1}{384}}} \Big[ -\,3\,H_2(H_1(H_2(\Phi_2))) -
3\,H_1(H_2(H_1(\Phi_1))) + 5\,H_1(H_2(H_2(\Phi_2))) +
5\,H_2(H_1(H_1(\Phi_1))) +4\,\Phi_1\,H_1(H_1(\Phi_2))+
\\
& \ \ \ \ \ \ \ \ \ \ \ \ \  + 4\,\Phi_2\,H_2(H_1(\Phi_2)) -
3\,\Phi_2\,H_1(H_1(\Phi_1)) -
3\,\Phi_1\,H_2(H_2(\Phi_2))-\,7\,\Phi_2\,H_1(H_2(\Phi_2)) -
7\,\Phi_1\,H_2(H_1(\Phi_1)) -
\\
& \ \ \ \ \ \ \ \ \ \ \ \ \  - 2\,H_1(\Phi_1)\,H_1(\Phi_2)
-2\,H_2(\Phi_2)\,H_2(\Phi_1) +4\,\Phi_1\Phi_2\,H_1(\Phi_1) +
4\,\Phi_1\Phi_2\,H_2(\Phi_2) \Big].
\endaligned
\end{equation*}
This geometry is equivalent to that of its model $\mathbb{ H}^3$ if
and only if its two essential curvatures $\kappa^{h_1t}_{i_1}$ and
$\kappa^{h_1t}_{i_2}$ vanish identically; equivalently, the two
explicit real functions $\mathbf{\Delta_1}$ and $\mathbf{\Delta_4}$
of only the three horizontal real variables $(x,y,u)$, with $z=x+iy,
w=u+iv$, vanish identically.
\end{Theorem}

Inspecting the method of construction of the fundamental vector
fields $H_1$ and $H_2$ in section 5 of \cite{Merker-Sabzevari-CEJM}
shows that they are in fact the real and imaginary parts of the
tangent vector field $2\overline{\mathcal L}$, introduced in this
paper. Moreover, checking the expressions of $T,\Phi_1,\Phi_2$ in
\cite{Merker-Sabzevari-CEJM}, enjoying the equalities:
\[
[H_1,H_2]=4T, \ \ \ \ [H_1,T]=\Phi_1\,T, \ \ \ \ [H_2,T]=\Phi_2\,T,
\]
specifies that we have:
\[\footnotesize\aligned
\mathcal L=\frac{1}{2}\,H_1-\frac{i}{2}\,H_2, \ \ \ \
\overline{\mathcal L}&=\frac{1}{2}\,H_1+\frac{i}{2}\,H_2, \ \ \ \
\mathcal T=-4\,T,
\\
P&=\frac{1}{2}\,\Phi_1-\frac{i}{2}\,\Phi_2.
\endaligned
\]
Now, putting the above complex expressions of $\mathcal
L,\overline{\mathcal L},\mathcal T,P$, into the single complex
essential invariant $\frak J$ of the equivalence problem of real
hypersurfaces $M^3\subset\mathbb C^2$ and comparing them carefully to
the above real expressions of the essential invariants
$\mathbf{\Delta}_1$ and $\mathbf{\Delta}_4$ of their Cartan
geometries surprisingly reveals that:

\begin{Theorem}
\label{Th-Geometry} The following relation holds between essential
invariants of the equivalence problem and Cartan geometry of the
Levi-nondegenerate $\mathcal C^6$-smooth real hypersurfaces
$M^3\subset\mathbb C^2$:
\[
\frak I=\frac{4}{{\sf c}\overline{\sf
c}^3}\big(\mathbf{\Delta}_1+i\,\mathbf{\Delta}_4\big).
\]
\end{Theorem}

This result shows that how much {\em explicitly} the two concepts of
equivalence problem and of Cartan geometry match up.

\end{document}